\documentclass[11pt]{amsart}
\usepackage{amsmath,amssymb,mathrsfs}
\newtheorem{theorem}{Theorem}[section]
\newtheorem{corollary}[theorem]{Corollary}
\newtheorem{lemma}[theorem]{Lemma}
\newtheorem{proposition}[theorem]{Proposition}

\begin{document}

\title{Ricci flat K\"ahler metrics with edge singularities}
\author{Simon Brendle}
\begin{abstract}
We construct Ricci flat K\"ahler metrics with cone singularities along a complex hypersurface. This construction is inspired in part by R.~Mazzeo's program in the case of negative Einstein constant, and uses the linear theory developed recently by S.~Donaldson.
\end{abstract}
\maketitle

\section{Introduction}

Let $(M,\omega_0)$ be a compact K\"ahler manifold of complex dimension $n$, and let $\Sigma$ be a smooth complex hypersurface in $M$. Moreover, let $\Lambda$ denote the holomorphic line bundle associated to $\Sigma$, and let $h$ be a bundle metric on $\Lambda$. We are interested in K\"ahler metrics on $M \setminus \Sigma$ which have edge singularities along $\Sigma$. To construct an example of such a metric, let us fix a smooth K\"ahler metric $\omega_0$ on $M$; a holomorphic section $s$ of $\Lambda$; and a real number $\beta \in (0,\frac{1}{2})$. If $\lambda > 0$ is sufficiently small, then the $(1,1)$-form 
\[\omega = \omega_0 + \lambda \, \sqrt{-1} \, \partial \bar{\partial}(|s|_h^{2\beta})\] 
defines a K\"ahler metric on $M \setminus \Sigma$ which has edge singularities along $\Sigma$. The number $\beta$ has a geometric interpretation in terms of the cone angle.

Our main result is an existence result for Ricci flat K\"ahler metrics with edge singularities. 

\begin{theorem}
\label{main.theorem}
Let $(M,\omega_0)$ be a compact K\"ahler manifold of complex dimension $n$. Moreover, let $\Sigma$ be a smooth complex hypersurface in $M$; let $\Lambda$ denote the holomorphic line bundle associated to $\Sigma$; and let $\beta \in (0,\frac{1}{2})$. Moreover, we assume that 
\[c_1(M) = (1-\beta) \, c_1(\Lambda) \in H^{1,1}(M,\mathbb{R}).\] 
Then there exists a K\"ahler metric $\hat{\omega} = \omega + \sqrt{-1} \, \partial\bar{\partial} u$ on $M \setminus \Sigma$ with the following properties: 
\begin{itemize} 
\item The metric $\hat{\omega}$ is Ricci flat with bounded curvature.
\item The metric $\hat{\omega}$ is uniformly equivalent to the background metric $\omega$, i.e. $a_1 \, \omega \leq \hat{\omega} \leq a_2 \, \omega$ for uniform constants $a_1,a_2 > 0$.
\end{itemize}
\end{theorem}

R.~Mazzeo \cite{Mazzeo2} has proposed a program for constructing K\"ahler-Einstein metrics with edge singularities in the case of negative Einstein constant. In dimension $2$, this problem has been studied by M.~Troyanov \cite{Troyanov}. Mazzeo's approach relies on the edge calculus developed in \cite{Mazzeo1}. The edge calculus has proved to be a very useful tool in the study of various geometric problems. We also note that G.~Tian and S.T~Yau have constructed complete Ricci flat K\"ahler metrics on the complement of a divisor; see \cite{Tian-Yau1}, \cite{Tian-Yau2} for details.

We now sketch the main steps involved in the proof of Theorem \ref{main.theorem}. In Section \ref{linear.theory}, we review some key results established in Donaldson's paper \cite{Donaldson}. Of particular importance is the Schauder estimate proved in Donaldson's paper. This estimate plays a key role in the argument.

In Section \ref{background.metric}, we describe the construction of the background metric $\omega$, and describe its basic properties. In particular, we show that $\text{\rm Ric}_\omega = \sqrt{-1} \, \partial\bar{\partial} F$ for some H\"older continuous function $F$. Moreover, we prove that $\omega$ has bounded curvature. Finally, we show that the covariant derivative of the curvature tensor of $\omega$ is bounded by $|DR|_g \leq C \, |\zeta|^{\varepsilon-\beta}$ for some $\varepsilon > 0$.

It is well known that the metric $\hat{\omega} = \omega + \sqrt{-1} \, \partial\bar{\partial} u$ is Ricci flat if and only if the function $u$ satisfies the complex Monge-Amp\`ere equation 
\[(\omega + \sqrt{-1} \, \partial\bar{\partial} u)^n = e^{F-c} \, \omega^n\] 
for some constant $c$. Following Aubin \cite{Aubin1} and Yau \cite{Yau2}, we solve this equation using the continuity method. To that end, we consider the equation 
\begin{equation} 
\tag{$\star_t$}
(\omega + \sqrt{-1} \, \partial\bar{\partial} u)^n = e^{tF-c} \, \omega^n, 
\end{equation}
where $t \in [0,1]$ is a parameter. 

In Section \ref{C0.estimate}, we establish a uniform $C^0$-estimate for solutions of $(\star_t)$. More precisely, we show that $\sup_M u - \inf_M u \leq C$ for any solution $u \in \mathcal{C}^{2,\alpha,\beta}$ of $(\star_t)$.

In Section \ref{C2.estimate}, we establish a Laplacian estimate for solutions of $(\star_t)$. In particular, this estimate implies that $a_1 \, \omega \leq \hat{\omega} \leq a_2 \, \omega$ whenever $u \in \mathcal{C}^{2,\alpha,\beta}$ is a solution of $(\star_t)$. In order to prove the Laplacian estimate, we adapt the arguments used in Yau's proof of the Calabi conjecture (cf. \cite{Yau1}, \cite{Yau2}). Our situation is more subtle in that we are dealing with a singular background metric $\omega$. Furthermore, the application of the maximum principle is not entirely straightforward, as the maximum may be attained on the singular set. In order to make the maximum principle work, we use a trick due to T.~Jeffres \cite{Jeffres1}, \cite{Jeffres2}. 

In Section \ref{C3.estimate}, we obtain an a-priori estimate for the covariant derivative of $\partial\bar{\partial} u$. This estimate is again proved using the maximum principle. In order to apply the maximum principle, we need to analyze the asymptotic behavior of the third order covariant derivatives of $u$ near $\Sigma$. 

Finally, in Section \ref{proof.of.main.theorem}, we show that the equation $(\star_1)$ has a solution $u \in \mathcal{C}^{2,\alpha,\beta}$. Moreover, we show that the associated K\"ahler metric $\hat{\omega} = \omega + \sqrt{-1} \, \partial\bar{\partial} u$ has bounded curvature, thereby completing the proof of Theorem \ref{main.theorem}.

\section{Donaldson's work on the linear theory}

\label{linear.theory}

In this section, we collect some results proved in Donaldson's paper \cite{Donaldson}, which play a major role in our argument. These results are all due to Donaldson, and are described here for the convenience of the reader. To begin with, we consider the K\"ahler metric 
\begin{align*} 
\Omega 
&= \sqrt{-1} \, \partial\bar{\partial}(|z_1|^2 + \hdots + |z_{n-1}|^2 + |\zeta|^{2\beta}) \\ 
&= \sum_{k=1}^{n-1} \sqrt{-1} \, dz_k \wedge d\bar{z}_k + \beta^2 \, \sqrt{-1} \, |\zeta|^{2\beta-2} \, d\zeta \wedge d\bar{\zeta} 
\end{align*}
on $\mathbb{C}^n$. The Laplacian with respect to that metric is given by 
\[\Delta_\Omega v = \sum_{k=1}^n \frac{\partial^2 v}{\partial z_k \, \partial \bar{z}_k} + \frac{1}{\beta^2} \, |\zeta|^{2-2\beta} \, \frac{\partial^2 v}{\partial \zeta \, \partial \bar{\zeta}}.\] 
The goal is to solve the equation $\Delta_\Omega v = f$. To that end, we write 
\[v(z_1,\hdots,z_{n-1},\zeta) = \tilde{v}(z_1,\hdots,z_{n-1},|\zeta|^{\beta-1} \, \zeta)\] 
and 
\[f(z_1,\hdots,z_{n-1},\zeta) = \tilde{f}(z_1,\hdots,z_{n-1},|\zeta|^{\beta-1} \, \zeta).\] 
The equation $\Delta_\Omega v = f$ is then equivalent to the equation 
\begin{equation} 
\label{cone.laplacian}
\sum_{k=1}^{n-1} \frac{\partial^2 \tilde{v}}{\partial z_k \, \partial \bar{z}_k} + \frac{1}{4} \, \Big ( \frac{\partial^2 \tilde{v}}{\partial r^2} + \frac{1}{r} \, \frac{\partial \tilde{v}}{\partial r} + \frac{1}{\beta^2r^2} \, \frac{\partial^2 \tilde{v}}{\partial \theta^2} \Big ) = \tilde{f}.
\end{equation}
Here, $|\zeta|^{\beta-1} \, \zeta = \xi = re^{i\theta}$. 

\begin{proposition}[S.~Donaldson \cite{Donaldson}]
Let $\tilde{f}$ be a function of class $C^\alpha$ with compact support. Then the equation (\ref{cone.laplacian}) has a unique weak solution $\tilde{v} \in L^{\frac{2n}{n-1}}(\mathbb{C}^n)$. The function $\tilde{v}$ is of class $C^{1,\alpha}$ and satisfies 
\begin{equation} 
\label{vanishing}
\frac{\partial \tilde{v}}{\partial \xi}(z_1,\hdots,z_{n-1},0) = \frac{\partial \tilde{v}}{\partial \bar{\xi}}(z_1,\hdots,z_{n-1},0) = 0. 
\end{equation}
Moreover, we have the estimate 
\begin{align} 
\label{schauder.estimate}
&\sum_{1 \leq k,l \leq n-1} \Big [ \frac{\partial^2 \tilde{v}}{\partial z_k \, \partial z_l} \Big ]_{C^\alpha} + \sum_{1 \leq k,l \leq n-1} \Big [ \frac{\partial^2 \tilde{v}}{\partial z_k \, \partial \bar{z}_l} \Big ]_{C^\alpha} \notag \\ 
&+ \sum_{k=1}^{n-1} \Big [ \frac{\partial^2 \tilde{v}}{\partial z_k \, \partial \xi} \Big ]_{C^\alpha} + \sum_{k=1}^{n-1} \Big [ \frac{\partial^2 \tilde{v}}{\partial z_k \, \partial \bar{\xi}} \Big ]_{C^\alpha} \leq C \, [\tilde{f}]_{C^\alpha}. 
\end{align} 
\end{proposition}

The identity (\ref{vanishing}) follows from the polyhomogeneous expansion for the Green's function derived in Donaldson's paper. By differentiating identity (\ref{vanishing}) in tangential direction, Donaldson showed that 
\begin{equation} 
\label{vanishing.2}
\frac{\partial^2 \tilde{v}}{\partial z_k \, \partial \xi}(z_1,\hdots,z_{n-1},0) = \frac{\partial^2 \tilde{v}}{\partial z_k \, \partial \bar{\xi}}(z_1,\hdots,z_{n-1},0) = 0 
\end{equation}
for each $k \in \{1,\hdots,n-1\}$.

Let us recall some notation from Donaldson's paper \cite{Donaldson}. Let $\mathcal{H}$ be the space of all functions $f$ of the form $f(z_1,\hdots,z_{n-1},\zeta) = \tilde{f}(z_1,\hdots,z_{n-1},|\zeta|^{\beta-1} \, \zeta)$, where $\tilde{f} \in W^{1,2}$. As in \cite{Donaldson}, we denote by $\mathcal{C}^{,\alpha,\beta}$ the space of all functions $f$ of the form $f(z_1,\hdots,z_{n-1},\zeta) = \tilde{f}(z_1,\hdots,z_{n-1},|\zeta|^{\beta-1} \, \zeta)$ where $\tilde{f} \in C^\alpha$. Moreover, let us denote by $\mathcal{C}_0^{,\alpha,\beta}$ the space of all functions $f \in \mathcal{C}^{,\alpha,\beta}$ such that $f(z_1,\hdots,z_{n-1},0) = 0$. A $(1,0)$-form $\tau$ is said to be of class $\mathcal{C}^{,\alpha,\beta}$ if 
\begin{align*} 
&\tau \Big ( \frac{\partial}{\partial z_k} \Big ) \in \mathcal{C}^{,\alpha,\beta}, \\ 
&|\zeta|^{1-\beta} \, \tau \Big ( \frac{\partial}{\partial \zeta} \Big ) \in \mathcal{C}_0^{,\alpha,\beta}. 
\end{align*}
Moreover, a $(1,1)$-form $\sigma$ is said to be of class $\mathcal{C}^{,\alpha,\beta}$ if 
\begin{align*} 
&\sigma \Big ( \frac{\partial}{\partial z_k},\frac{\partial}{\partial \bar{z}_l} \Big ) \in \mathcal{C}^{,\alpha,\beta}, \\ 
&|\zeta|^{1-\beta} \, \sigma \Big ( \frac{\partial}{\partial z_k},\frac{\partial}{\partial \bar{\zeta}} \Big ) \in \mathcal{C}_0^{,\alpha,\beta}, \\ 
&|\zeta|^{1-\beta} \, \sigma \Big ( \frac{\partial}{\partial \zeta},\frac{\partial}{\partial \bar{z}_l} \Big ) \in \mathcal{C}_0^{,\alpha,\beta}, \\ 
&|\zeta|^{2-2\beta} \, \sigma \Big ( \frac{\partial}{\partial \zeta},\frac{\partial}{\partial \bar{\zeta}} \Big ) \in \mathcal{C}^{,\alpha,\beta}. 
\end{align*}
Finally, let 
\[\mathcal{C}^{2,\alpha,\beta} = \{v \in C^2(M \setminus \Sigma): \text{\rm $v,\partial v,\partial\bar{\partial} v$ are of class $\mathcal{C}^{,\alpha,\beta}$}\}\] 
(cf. \cite{Donaldson}, p.~16). 

Consider now a function $f \in \mathcal{C}^{,\alpha,\beta}$ which is supported in the ball $B_6(0)$. In view of Donaldson's results, the equation $\Delta_\Omega v = f$ admits a weak solution $v \in \mathcal{H}$ satisfying $\int_{\mathbb{C}^n} |v|^{\frac{2n}{n-1}} \, \Omega^n < \infty$. The Schauder estimate (\ref{schauder.estimate}) implies 
\[\sum_{1 \leq k,l \leq n-1} \Big \| \frac{\partial^2 v}{\partial z_k \, \partial z_l} \Big \|_{\mathcal{C}^{,\alpha,\beta}} + \sum_{1 \leq k,l \leq n-1} \Big \| \frac{\partial^2 v}{\partial z_k \, \partial \bar{z}_l} \Big \|_{\mathcal{C}^{,\alpha,\beta}} \leq C \, \|f\|_{\mathcal{C}^{,\alpha,\beta}}.\] 
Using the equation $\Delta_\Omega v = f$, Donaldson obtained the estimate 
\[\Big \| |\zeta|^{2\beta-2} \, \frac{\partial^2 v}{\partial \zeta \, \partial \bar{\zeta}} \Big \|_{\mathcal{C}^{,\alpha,\beta}} \leq C \, \|f\|_{\mathcal{C}^{,\alpha,\beta}}.\] 
We next observe that 
\begin{align*} 
|\zeta|^{1-\beta} \, \frac{\partial^2 v}{\partial z_k \, \partial \bar{\zeta}}(z_1,\hdots,z_{n-1},\zeta) 
&= \frac{\beta-1}{2} \, \frac{\zeta^2}{|\zeta|^2} \, \frac{\partial^2 \tilde{v}}{\partial z_k \, \partial \xi}(z_1,\hdots,z_{n-1},|\zeta|^{\beta-1} \, \zeta) \\ 
&+ \frac{\beta+1}{2} \, \frac{\partial^2 \tilde{v}}{\partial z_k \, \partial \bar{\xi}}(z_1,\hdots,z_{n-1},|\zeta|^{\beta-1} \, \zeta) 
\end{align*} 
for $k \in \{1,\hdots,n-1\}$. In view of (\ref{vanishing.2}), the functions 
\[(z_1,\hdots,z_{n-1},\zeta) \mapsto \frac{\partial^2 \tilde{v}}{\partial z_k \, \partial \xi}(z_1,\hdots,z_{n-1},|\zeta|^{\beta-1} \, \zeta)\] 
and 
\[(z_1,\hdots,z_{n-1},\zeta) \mapsto \frac{\partial^2 \tilde{v}}{\partial z_k \, \partial \bar{\xi}}(z_1,\hdots,z_{n-1},|\zeta|^{\beta-1} \, \zeta)\] 
belong to the space $\mathcal{C}_0^{,\alpha,\beta}$. Consequently, the function $|\zeta|^{1-\beta} \, \frac{\partial^2 v}{\partial z_k \, \partial \bar{\zeta}}$ belongs to the space $\mathcal{C}_0^{,\alpha,\beta}$, and one has the estimate 
\[\Big \| |\zeta|^{1-\beta} \, \frac{\partial^2 v}{\partial z_k \, \partial \bar{\zeta}} \Big \|_{\mathcal{C}^{,\alpha,\beta}} \leq C \, \|f\|_{\mathcal{C}^{,\alpha,\beta}}.\] 
Donaldson's Schauder estimate can be summarized as follows: 

\begin{proposition}[S.~Donaldson \cite{Donaldson}]
\label{model.problem}
Let $f \in \mathcal{C}^{,\alpha,\beta}$ be a function supported in the ball $B_6(0)$. Then there exists a unique weak solution $v \in \mathcal{H}$ of the equation $\Delta_\Omega v = f$ such that $\int_{\mathbb{C}^n} |v|^{\frac{2n}{n-1}} \, \Omega^n < \infty$. Moreover, $v$ is of class $\mathcal{C}^{2,\alpha,\beta}$ and we have the estimate 
\[\|v\|_{\mathcal{C}^{2,\alpha,\beta}} \leq C \, \|f\|_{\mathcal{C}^{,\alpha,\beta}}.\] 
\end{proposition}

It is shown in Section 4.2 of Donaldson's paper that the Schauder estimate in Proposition \ref{model.problem} carries over to the variable coefficient setting. In the following, $\omega$ will denote a K\"ahler metric on the ball $B_6(0)$ with the property that $\omega$ is of class $\mathcal{C}^{,\alpha,\beta}$ and $\|\omega - \Omega\|_{\mathcal{C}^{,\alpha,\beta}} \leq \varepsilon_0$. Here, $\varepsilon_0$ denotes a small positive constant.

\begin{proposition}[S.~Donaldson \cite{Donaldson}]
\label{interior.estimate}
Let $\omega$ be a K\"ahler metric on the ball $B_6(0)$ such that $\omega$ is of class $\mathcal{C}^{,\alpha,\beta}$ and $\|\omega - \Omega\|_{\mathcal{C}^{,\alpha,\beta}} \leq \varepsilon_0$. Moreover, let $v$ be a function of class $\mathcal{C}^{2,\alpha,\beta}$ defined on the ball $B_2(0)$. If $\varepsilon_0$ is sufficiently small, then we have the estimate 
\[\|v\|_{\mathcal{C}^{2,\alpha,\beta}(B_1(0))} \leq C \, \sup_{B_2(0)} |v| + C \, \|\Delta_\omega v\|_{\mathcal{C}^{,\alpha,\beta}(B_2(0))}.\] 
\end{proposition} 

Proposition \ref{interior.estimate} is a consequence of Proposition \ref{model.problem} (compare \cite{Simon}). Note that Proposition \ref{interior.estimate} presupposes that the function $v$ is of class $\mathcal{C}^{2,\alpha,\beta}$. This assumption can be removed as follows: 

\begin{proposition}[S.~Donaldson \cite{Donaldson}]
\label{regularity.1}
Let $\omega$ be a K\"ahler metric on the ball $B_6(0)$ which is of class $\mathcal{C}^{,\alpha,\beta}$ and satisfies $\|\omega - \Omega\|_{\mathcal{C}^{,\alpha,\beta}} \leq \varepsilon_0$. Suppose that $\varphi$ is a function of class $\mathcal{C}^{,\alpha,\beta}$ defined on the ball $B_6(0)$, and $v \in \mathcal{H}$ is a weak solution of the equation $\Delta_\omega v = \varphi$. If $\varepsilon_0$ is sufficiently small, then the restriction $v|_{B_1(0)}$ is of class $\mathcal{C}^{2,\alpha,\beta}$. 
\end{proposition}

We define 
\[\mathcal{X} = \{f \in \mathcal{C}^{,\alpha,\beta}: \text{\rm $f$ vanishes on $\mathbb{C}^n \setminus B_6(0)$}\}.\] 
Given any function $f \in \mathcal{X}$, we denote by $v = Gf$ the unique weak solution of the equation $\Delta_\Omega v = f$. By Proposition \ref{model.problem}, this defines a bounded linear operator $G: \mathcal{X} \to \mathcal{C}^{2,\alpha,\beta}$.

Let us fix a smooth cutoff function $\chi$ such that $\chi=1$ in $B_5(0)$ and $\chi=0$ on $\mathbb{C}^n \setminus B_6(0)$. We define a differential operator $L$ by 
\[Lv = \sqrt{-1} \, \frac{\partial((1-\chi) \, \bar{\partial} v) \wedge \Omega^{n-1} + \partial(\chi \, \bar{\partial} v) \wedge \omega^{n-1}}{(1-\chi) \, \Omega^n + \chi \, \omega^n}.\] 
Clearly, $Lv = \Delta_\Omega v$ on $\mathbb{C}^n \setminus B_6(0)$. Consequently, the operator $LG: \mathcal{X} \to \mathcal{C}^{,\alpha,\beta}$ maps $\mathcal{X}$ into itself. Moreover, we have $\|LGf - f\|_{\mathcal{C}^{,\alpha,\beta}} \leq C\varepsilon_0 \, \|f\|_{\mathcal{C}^{,\alpha,\beta}}$. Hence, if $\varepsilon_0$ is small enough, then the operator $LG: \mathcal{X} \to \mathcal{X}$ is invertible. Let $T: \mathcal{X} \to \mathcal{X}$ denote the inverse of $LG$. 

We now sketch the proof of Proposition \ref{regularity.1}. Let $\eta$ be a smooth cutoff function satisfying $\eta=1$ in $B_4(0)$ and $\eta=0$ in $\mathbb{C}^n \setminus B_5(0)$, and let $\psi = \eta\varphi + 2 \, \langle d\eta,dv \rangle_\omega + v \, \Delta_\omega \eta$. Then $\int_{B_6(0)} |\psi|^2 \, \Omega^n < \infty$ and the restriction $\psi|_{B_4(0)}$ is of class $\mathcal{C}^{,\alpha,\beta}$. We can find a sequence of functions $\psi_k \in \mathcal{X}$ such that $\psi_k|_{B_3(0)} = \psi|_{B_3(0)}$ and 
\begin{equation} 
\label{approx.1}
\int_{B_6(0)} |\psi - \psi_k|^2 \, \Omega^n \to 0 
\end{equation}
as $k \to \infty$. For each $k$, the function $v_k = GT\psi_k$ is of class $\mathcal{C}^{2,\alpha,\beta}$, and we have $Lv_k = \psi_k$ away from the set $\{\zeta=0\}$. On the other hand, the function $\eta v \in \mathcal{H}$ is a weak solution of the equation $L(\eta v) = \psi$. Consequently, the function $\eta v-v_k \in \mathcal{H}$ is a weak solution of the equation $L(\eta v - v_k) = \psi-\psi_k$. Using (\ref{approx.1}), we obtain 
\begin{equation} 
\label{approx.2}
\int_{\mathbb{C}^n} |\eta v - v_k|^{\frac{2n}{n-1}} \, \Omega^n \to 0 
\end{equation}
as $k \to \infty$. 

Note that the pull-back of the K\"ahler metric $\omega$ under the map $\Phi: (z_1,\hdots,z_{n-1},\xi) \mapsto (z_1,\hdots,z_{n-1},|\xi|^{\frac{1}{\beta}-1} \, \xi)$ is uniformly equivalent to the Euclidean metric. Hence, it follows from Theorem 8.17 in \cite{Gilbarg-Trudinger} that 
\[\sup_{B_2(0)} |v_k| \leq C \, \bigg ( \int_{B_3(0)} |v_k|^{\frac{2n}{n-1}} \, \Omega^n \bigg )^{\frac{n-1}{2n}} + C \, \sup_{B_3(0)} |\Delta_\omega v_k|.\] 
Moreover, since $v_k \in \mathcal{C}^{2,\alpha,\beta}$, we have 
\[\|v_k\|_{\mathcal{C}^{2,\alpha,\beta}(B_1(0))} \leq C \, \sup_{B_2(0)} |v_k| + C \, \|\Delta_\omega v_k\|_{\mathcal{C}^{,\alpha,\beta}(B_2(0))}\] 
by Proposition \ref{interior.estimate}. On the other hand, we have $\Delta_\omega v_k = Lv_k = \psi_k = \psi$ in $B_3(0)$. Putting these facts together, we obtain 
\[\|v_k\|_{\mathcal{C}^{2,\alpha,\beta}(B_1(0))} \leq C \, \bigg ( \int_{B_3(0)} |v_k|^{\frac{2n}{n-1}} \, \Omega^n \bigg )^{\frac{n-1}{2n}} + C \, \|\psi\|_{\mathcal{C}^{,\alpha,\beta}(B_3(0))}\] 
for some uniform constant $C$. Using (\ref{approx.2}), we conclude that $v|_{B_1(0)} \in \mathcal{C}^{2,\alpha,\beta}$, as claimed. \\

By freezing coefficients, Donaldson obtained the following regularity result: 

\begin{theorem}[S.~Donaldson \cite{Donaldson}]
\label{regularity.2}
Let $\omega$ be a K\"ahler metric on the ball $B_6(0)$ which is of class $\mathcal{C}^{,\alpha,\beta}$ and satisfies $a_1 \, \Omega \leq \omega \leq a_2 \, \Omega$ for suitable constants $a_1,a_2>0$. Suppose that $\varphi$ is a function of class $\mathcal{C}^{,\alpha,\beta}$ defined on the ball $B_6(0)$, and $v \in \mathcal{H}$ is a weak solution of the equation $\Delta_\omega v = \varphi$. Then the restriction $v|_{B_1(0)}$ is of class $\mathcal{C}^{2,\alpha,\beta}$.
\end{theorem}

The following is a special case of Theorem \ref{regularity.2}.

\begin{corollary}[S.~Donaldson \cite{Donaldson}] 
\label{regularity.3}
Let $\omega$ be a K\"ahler metric on the ball $B_6(0)$ which is of class $\mathcal{C}^{,\alpha,\beta}$ and satisfies $a_1 \, \Omega \leq \omega \leq a_2 \, \Omega$ for suitable constants $a_1,a_2>0$. Suppose that $v$ and $\varphi$ are functions of class $\mathcal{C}^{,\alpha,\beta}$ defined on the ball $B_6(0)$ such that $\Delta_\omega v = \varphi$ away from the set $\{\zeta=0\}$. Then the restriction $v|_{B_1(0)}$ is of class $\mathcal{C}^{2,\alpha,\beta}$.
\end{corollary}


\section{The background edge metric}

\label{background.metric}

In the following, we fix a K\"ahler manifold $(M,\omega_0)$, a smooth complex hypersurface $\Sigma \subset M$, and a real number $\beta \in (0,\frac{1}{2})$. Let $\Lambda$ be the holomorphic line bundle associated with $\Sigma$, and let $s$ be a holomorphic section of $\Lambda$ such that $\Sigma = \{s=0\}$. 

The $(1,1)$-form $-\sqrt{-1} \, \partial\bar{\partial} \log(|s|_h^2)$ is smooth and represents the cohomology class $c_1(\Lambda) \in H^{1,1}(M,\mathbb{R})$. Consequently, the cohomology class $c_1(M) - (1-\beta) \, c_1(\Lambda) \in H^{1,1}(M,\mathbb{R})$ can be represented by the $(1,1)$-form 
\[\text{\rm Ric}_{\omega_0} + (1-\beta) \, \sqrt{-1} \, \partial\bar{\partial} \log(|s|_h^2).\] 
Since $c_1(M) = (1-\beta) \, c_1(\Lambda)$, there exists a smooth function $F_0$ such that 
\[\text{\rm Ric}_{\omega_0} + (1-\beta) \, \sqrt{-1} \, \partial\bar{\partial} \log(|s|_h^2) = \sqrt{-1} \, \partial\bar{\partial} F_0.\] 
We next define 
\[\omega = \omega_0 + \lambda \, \sqrt{-1} \, \partial\bar{\partial} (|s|_h^{2\beta}).\] 
It is easy to see that 
\[\omega \geq \omega_0 + \lambda\beta \, |s|_h^{2\beta} \, \sqrt{-1} \, \partial\bar{\partial} \log(|s|_h^2).\] 
The term $\partial\bar{\partial} \log(|s|_h^2)$ defines a smooth $(1,1)$-form on $M$. Hence, if we choose $\lambda > 0$ sufficiently small, then $\omega$ is a K\"ahler metric on $M \setminus \Sigma$. The Ricci curvature of $\omega$ is given by 
\begin{align*} 
\text{\rm Ric}_\omega 
&= \text{\rm Ric}_{\omega_0} - \sqrt{-1} \, \partial \bar{\partial} \log \frac{\omega^n}{\omega_0^n} \\ 
&= \sqrt{-1} \, \partial\bar{\partial} F_0 - (1-\beta) \, \sqrt{-1} \, \partial\bar{\partial} \log(|s|_h^2) - \sqrt{-1} \, \partial \bar{\partial} \log \frac{\omega^n}{\omega_0^n} \\ 
&= \sqrt{-1} \, \partial\bar{\partial} F,
\end{align*}
where 
\[F = F_0 - \log \frac{|s|_h^{2-2\beta} \, \omega^n}{\omega_0^n}.\] 
Clearly, $F$ is a smooth function away from $\Sigma$. We next analyze the behavior of $F$ near $\Sigma$. To that end, we fix a point $p \in \Sigma$ and introduce complex coordinates $(z_1,\hdots,z_{n-1},\zeta)$ near $p$ so that $\Sigma = \{\zeta=0\}$. Moreover, we write $|s|_h^{2\beta} = \rho \, |\zeta|^{2\beta}$, where $\rho$ is a smooth positive function. Then 
\begin{equation} 
\label{omega}
\omega = \omega_0 + \lambda \, \sqrt{-1} \, \partial\bar{\partial} (\rho \, |\zeta|^{2\beta}), 
\end{equation} 
hence 
\begin{align*} 
\omega 
&= \omega_0 + \lambda \, |\zeta|^{2\beta} \, \sqrt{-1} \, \partial\bar{\partial} \rho \\ 
&+ \lambda\beta \, |\zeta|^{2\beta-2} \, \zeta \, \sqrt{-1} \, \partial \rho \wedge d\bar{\zeta} \\ 
&+ \lambda\beta \, |\zeta|^{2\beta-2} \, \bar{\zeta} \, \sqrt{-1} \, d\zeta \wedge \bar{\partial} \rho \\ 
&+ \lambda\beta^2\rho \, \sqrt{-1} \, |\zeta|^{2\beta-2} \, d\zeta \wedge d\bar{\zeta}. 
\end{align*}
This implies 
\[\frac{\omega^n}{\omega_0^n} = |\zeta|^{2\beta-2} \sum_{m=0}^{n-1} a_m \, |\zeta|^{2\beta m} + \sum_{m=0}^n b_m \, |\zeta|^{2\beta m},\] 
where $a_m,b_m$ are smooth functions and $a_0$ is positive. Consequently, 
\begin{align*} 
\frac{|s|_h^{2-2\beta} \, \omega^n}{\omega_0^n} 
&= \rho^{\frac{1}{\beta}-1} \, \frac{|\zeta|^{2-2\beta} \, \omega^n}{\omega_0^n} \\ 
&= \rho^{\frac{1}{\beta}-1} \, \bigg ( \sum_{m=0}^{n-1} a_m \, |\zeta|^{2\beta m} + |\zeta|^{2-2\beta} \sum_{m=0}^n b_m \, |\zeta|^{2\beta m} \bigg ), 
\end{align*}
hence 
\begin{align} 
\label{formula.F}
F &= F_0 - \Big ( \frac{1}{\beta}-1 \Big ) \, \log \rho \notag \\ 
&- \log \bigg ( \sum_{m=0}^{n-1} a_m \, |\zeta|^{2\beta m} + |\zeta|^{2-2\beta} \sum_{m=0}^n b_m \, |\zeta|^{2\beta m} \bigg ). 
\end{align}
Since $a_m$, $b_m$, and $\rho$ are smooth functions, we can draw the following conclusions:

\begin{proposition}
\label{properties.F}
The function $F$ is H\"older continuous. Moreover, we have 
\[|\partial\bar{\partial} F|_g \leq C + C \, |\zeta|^{2-4\beta}.\] 
In particular, if $\beta \in (0,\frac{1}{2})$, then $|\partial\bar{\partial} F|_g \leq C$. 
\end{proposition}

In the next step, we estimate the Riemann curvature tensor of $g$. As above, we may work in local complex coordinates around a point $p \in \Sigma$. We will consider the pull-back of $\omega$ under the holomorphic map 
\[\Psi: (z_1,\hdots,z_{n-1},z_n) \mapsto (z_1,\hdots,z_{n-1},z_n^{\frac{1}{\beta}}).\] 
Note that $\Psi$ is defined only locally. Using the identity (\ref{omega}), we obtain 
\begin{equation} 
\Psi^*\omega = \Psi^*\omega_0 + \lambda \, \sqrt{-1} \, \partial\bar{\partial} \big ( (\rho \circ \Psi) \, |z_n|^2 \big ), 
\end{equation}
hence 
\begin{align*} 
\Psi^*\omega 
&= \Psi^*\omega_0 + \lambda \, |z_n|^2 \, \sqrt{-1} \, \partial\bar{\partial} (\rho \circ \Psi) \\ 
&+ \lambda \, z_n \, \sqrt{-1} \, \partial (\rho \circ \Psi) \wedge d\bar{z}_n \\ 
&+ \lambda \, \bar{z}_n \, \sqrt{-1} \, dz_n \wedge \bar{\partial} (\rho \circ \Psi) \\ 
&+ \lambda \, (\rho \circ \Psi) \, \sqrt{-1} \, dz_n \wedge d\bar{z}_n. 
\end{align*}
This implies 
\[\Psi^*g \Big ( \frac{\partial}{\partial z_k},\frac{\partial}{\partial \bar{z}_l} \Big ) = (g_{0,k\bar{l}} \circ \Psi) + \lambda \, |z_n|^2 \, \Big ( \frac{\partial^2 \rho}{\partial z_k \, \partial \bar{z}_l} \circ \Psi \Big )\] 
and 
\begin{align*} 
\Psi^*g \Big ( \frac{\partial}{\partial z_n},\frac{\partial}{\partial \bar{z}_l} \Big ) 
&= \frac{1}{\beta} \, z_n^{\frac{1}{\beta}-1} \, (g_{0,\zeta\bar{l}} \circ \Psi) + \frac{\lambda}{\beta} \, 
z_n^{\frac{1}{\beta}} \, \bar{z}_n \, \Big ( \frac{\partial^2 \rho}{\partial \zeta \, \partial \bar{z}_l} \circ \Psi \Big ) \\ 
&+ \lambda \, \bar{z}_n \, \Big ( \frac{\partial \rho}{\partial \bar{z}_l} \circ \Psi \Big ) 
\end{align*} 
for all $k,l \in \{1,\hdots,n-1\}$. Moreover, we have 
\begin{align*}
\Psi^*g \Big ( \frac{\partial}{\partial z_n},\frac{\partial}{\partial \bar{z}_n} \Big ) 
&= \frac{1}{\beta^2} \, |z_n|^{\frac{2}{\beta}-2} \, (g_{0,\zeta\bar{\zeta}} \circ \Psi) + \frac{\lambda}{\beta^2} \, |z_n|^{\frac{2}{\beta}} \, \Big ( \frac{\partial^2 \rho}{\partial \zeta \, \partial \bar{\zeta}} \circ \Psi \Big ) \\ 
&+ \frac{\lambda}{\beta} \, z_n^{\frac{1}{\beta}} \, \Big ( \frac{\partial \rho}{\partial \zeta} \circ \Psi \Big ) + \frac{\lambda}{\beta} \, \bar{z}_n^{\frac{1}{\beta}} \, \Big ( \frac{\partial \rho}{\partial \bar{\zeta}} \circ \Psi \Big ) + \lambda \, (\rho \circ \Psi).
\end{align*}
Using these identities, we obtain estimates for the first derivatives of $\Psi^* g$.

\begin{lemma}
\label{est.1}
We have 
\begin{align*} 
&\sum_{i,k,l=1}^{n-1} \Big | \frac{\partial}{\partial z_i} \Psi^*g \Big ( \frac{\partial}{\partial z_k},\frac{\partial}{\partial \bar{z}_l} \Big ) \Big | \leq C, \\ 
&\sum_{i,l=1}^{n-1} \Big | \frac{\partial}{\partial z_i} \Psi^*g \Big ( \frac{\partial}{\partial z_n},\frac{\partial}{\partial \bar{z}_l} \Big ) \Big | \leq C \, |z_n| + C \, |z_n|^{\frac{1}{\beta}-1}, \\ 
&\sum_{i=1}^{n-1} \Big | \frac{\partial}{\partial z_i} \Psi^*g \Big ( \frac{\partial}{\partial z_n},\frac{\partial}{\partial \bar{z}_n} \Big ) \Big | \leq C, \\ 
&\sum_{l=1}^{n-1} \Big | \frac{\partial}{\partial z_n} \Psi^*g \Big ( \frac{\partial}{\partial z_n},\frac{\partial}{\partial \bar{z}_l} \Big ) \Big | \leq C \, |z_n|^{\frac{1}{\beta}-2}, \\ 
&\Big | \frac{\partial}{\partial z_n} \Psi^*g \Big ( \frac{\partial}{\partial z_n},\frac{\partial}{\partial \bar{z}_n} \Big ) \Big | \leq C \, |z_n|^{\frac{1}{\beta}-1} + C \, |z_n|^{\frac{2}{\beta}-3}. 
\end{align*}
\end{lemma}

The second derivatives of $\Psi^* g$ can be estimated as follows:

\begin{lemma}
\label{est.2}
We have 
\begin{align*} 
&\sum_{i,j,k,l=1}^{n-1} \Big | \frac{\partial^2}{\partial z_i \, \partial \bar{z}_j} \Psi^*g \Big ( \frac{\partial}{\partial z_k},\frac{\partial}{\partial \bar{z}_l} \Big ) \Big | \leq C, \\
&\sum_{i,j,l=1}^{n-1} \Big | \frac{\partial^2}{\partial z_i \, \partial \bar{z}_j} \Psi^*g \Big ( \frac{\partial}{\partial z_n},\frac{\partial}{\partial \bar{z}_l} \Big ) \Big | \leq C \, |z_n| + C \, |z_n|^{\frac{1}{\beta}-1}, \\ 
&\sum_{i,j=1}^n \Big | \frac{\partial^2}{\partial z_i \, \partial \bar{z}_j} \Psi^*g \Big ( \frac{\partial}{\partial z_n},\frac{\partial}{\partial \bar{z}_n} \Big ) \Big | \leq C, \\
&\sum_{j,l=1}^n \Big | \frac{\partial^2}{\partial z_n \, \partial \bar{z}_j} \Psi^*g \Big ( \frac{\partial}{\partial z_n},\frac{\partial}{\partial \bar{z}_l} \Big ) \Big | \leq C \, |z_n|^{\frac{1}{\beta}-2}, \\
&\sum_{i=1}^n \Big | \frac{\partial^2}{\partial z_i \, \partial \bar{z}_n} \Psi^*g \Big ( \frac{\partial}{\partial z_n},\frac{\partial}{\partial \bar{z}_n} \Big ) \Big | \leq C \, |z_n|^{\frac{1}{\beta}-1} + C \, |z_n|^{\frac{2}{\beta}-3}, \\ 
&\Big | \frac{\partial^2}{\partial z_n \, \partial \bar{z}_n} \Psi^*g \Big ( \frac{\partial}{\partial z_n},\frac{\partial}{\partial \bar{z}_n} \Big ) \Big | \leq C \, |z_n|^{\frac{2}{\beta}-4}. 
\end{align*}
\end{lemma}

We also note that 
\[\Big | \frac{\partial^2}{\partial z_n \, \partial z_n} \Psi^*g \Big ( \frac{\partial}{\partial z_n},\frac{\partial}{\partial \bar{z}_l} \Big ) \Big | \leq C \, |z_n|^{\frac{1}{\beta}-3}.\] 
Note that this term does not enter into the formula for the curvature tensor.

After these preparations, we now derive a bound for the Riemann curvature tensor of $g$.

\begin{proposition}
\label{curvature}
The curvature tensor of $g$ can be estimated by 
\[|R|_g \leq C + C \, |z_n|^{\frac{1}{\beta}-2} + C \, |z_n|^{\frac{2}{\beta}-4}.\] 
In particular, if $\beta \in (0,\frac{1}{2})$, then the curvature tensor of $\omega$ is bounded. 
\end{proposition}

\textbf{Proof.} 
The curvature tensor of a K\"ahler manifold is given by 
\begin{equation} 
\label{kahler.curvature}
R_{\alpha\bar{\beta}\gamma\bar{\delta}} = -\partial_\alpha \partial_{\bar{\beta}} g_{\gamma\bar{\delta}} + \sum_{\mu,\nu=1}^n g^{\mu\bar{\nu}} \, \partial_\gamma g_{\alpha\bar{\nu}} \, \partial_{\bar{\delta}} \, g_{\mu\bar{\beta}} 
\end{equation}
(cf. \cite{Yau2}, page 344). Note that the metric $\Psi^*g$ is uniformly equivalent to the Euclidean metric. Hence, it follows from Lemma \ref{est.1} and Lemma \ref{est.2} that 
\begin{align*} 
|R|_g 
&\leq C \sum_{\alpha,\beta,\gamma,\delta}^n \Big | R \Big ( \frac{\partial}{\partial z_\alpha},\frac{\partial}{\partial \bar{z}_\beta},\frac{\partial}{\partial z_\gamma},\frac{\partial}{\partial \bar{z}_\delta} \Big ) \Big | \\ 
&\leq \sum_{\alpha,\beta,\gamma,\delta=1}^n \Big | \frac{\partial^2}{\partial z_\alpha \, \partial \bar{z}_\beta} \Psi^*g \Big ( \frac{\partial}{\partial z_\gamma},\frac{\partial}{\partial \bar{z}_\delta} \Big ) \Big | + C \, \sum_{\alpha,\gamma,\delta=1}^n \Big | \frac{\partial}{\partial z_\alpha} \Psi^*g \Big ( \frac{\partial}{\partial z_\gamma},\frac{\partial}{\partial \bar{z}_\delta} \Big ) \Big |^2 \\ 
&\leq C + C \, |z_n|^{\frac{1}{\beta}-2} + C \, |z_n|^{\frac{2}{\beta}-4}.
\end{align*}
This completes the proof. \\

A similar argument gives a bound for the covariant derivative of the curvature tensor.

\begin{proposition}
\label{derivative.of.curvature}
We have
\[|DR|_g \leq C + C \, |z_n|^{\frac{1}{\beta}-3} + C \, |z_n|^{\frac{2}{\beta}-5}.\] 
In particular, if $\beta \in (0,\frac{1}{2})$, then $|DR|_g \leq C \, |z_n|^{\frac{\varepsilon}{\beta}-1}$ for some $\varepsilon>0$.
\end{proposition}

\textbf{Proof.} 
It follows from Lemma \ref{est.1} and Lemma \ref{est.2} that 
\[\sum_{\alpha,\gamma,\delta=1}^n \Big | \frac{\partial}{\partial z_\alpha} \Psi^*g \Big ( \frac{\partial}{\partial z_\gamma},\frac{\partial}{\partial \bar{z}_\delta} \Big ) \Big | \leq C + C \, |z_n|^{\frac{1}{\beta}-2}\] 
and 
\[\sum_{\alpha,\beta,\gamma,\delta=1}^n \Big | \frac{\partial^2}{\partial z_\alpha \, \partial \bar{z}_\beta} \Psi^*g \Big ( \frac{\partial}{\partial z_\gamma},\frac{\partial}{\partial \bar{z}_\delta} \Big ) \Big | \leq C + C \, |z_n|^{\frac{1}{\beta}-2} + C \, |z_n|^{\frac{2}{\beta}-4}.\] 
A similar calculation gives 
\[\sum_{\alpha,\beta,\gamma,\delta=1}^n \Big | \frac{\partial^2}{\partial z_\alpha \, \partial z_\beta} \Psi^*g \Big ( \frac{\partial}{\partial z_\gamma},\frac{\partial}{\partial \bar{z}_\delta} \Big ) \Big | \leq C + C \, |z_n|^{\frac{1}{\beta}-3}\] 
and 
\[\sum_{\alpha,\beta,\gamma,\delta,\mu=1}^n \Big | \frac{\partial^3}{\partial z_\alpha \, \partial \bar{z}_\beta \, \partial z_\mu} \Psi^*g \Big ( \frac{\partial}{\partial z_\gamma},\frac{\partial}{\partial \bar{z}_\delta} \Big ) \Big | \leq C + C \, |z_n|^{\frac{1}{\beta}-3} + C \, |z_n|^{\frac{2}{\beta}-5}.\] 
Since the metric $\Psi^*g$ is uniformly equivalent to the Euclidean metric, we conclude that 
\begin{align*} 
|DR|_g &\leq C \, \sum_{\alpha,\beta,\gamma,\delta,\mu=1}^n \Big | \frac{\partial^3}{\partial z_\alpha \, \partial \bar{z}_\beta \, \partial z_\mu} \Psi^*g \Big ( \frac{\partial}{\partial z_\gamma},\frac{\partial}{\partial \bar{z}_\delta} \Big ) \Big | \\ 
&+ C \, \bigg ( \sum_{\alpha,\gamma,\delta=1}^n \Big | \frac{\partial}{\partial z_\alpha} \Psi^*g \Big ( \frac{\partial}{\partial z_\gamma},\frac{\partial}{\partial \bar{z}_\delta} \Big ) \Big | \bigg ) \, \bigg ( \sum_{\alpha,\beta,\gamma,\delta=1}^n \Big | \frac{\partial^2}{\partial z_\alpha \, \partial z_\beta} \Psi^*g \Big ( \frac{\partial}{\partial z_\gamma},\frac{\partial}{\partial \bar{z}_\delta} \Big ) \Big | \bigg ) \\ 
&+ C \, \bigg ( \sum_{\alpha,\gamma,\delta=1}^n \Big | \frac{\partial}{\partial z_\alpha} \Psi^*g \Big ( \frac{\partial}{\partial z_\gamma},\frac{\partial}{\partial \bar{z}_\delta} \Big ) \Big | \bigg ) \, \bigg ( \sum_{\alpha,\beta,\gamma,\delta=1}^n \Big | \frac{\partial^2}{\partial z_\alpha \, \partial \bar{z}_\beta} \Psi^*g \Big ( \frac{\partial}{\partial z_\gamma},\frac{\partial}{\partial \bar{z}_\delta} \Big ) \Big | \bigg ) \\ 
&+ C \, \bigg ( \sum_{\alpha,\gamma,\delta=1}^n \Big | \frac{\partial}{\partial z_\alpha} \Psi^*g \Big ( \frac{\partial}{\partial z_\gamma},\frac{\partial}{\partial \bar{z}_\delta} \Big ) \Big | \bigg )^3 \\ 
&\leq C + C \, |z_n|^{\frac{1}{\beta}-3} + C \, |z_n|^{\frac{2}{\beta}-5}.
\end{align*} 
This completes the proof. \\

We next describe a Sobolev inequality for the background edge manifold $(M \setminus \Sigma,g)$. This inequality will be used in the proof of the $C^0$ estimate.

\begin{proposition}
\label{sobolev}
Let $v$ be a smooth function on $M \setminus \Sigma$ satisfying $\sup_{M \setminus \Sigma} |v| < \infty$. Then 
\begin{equation} 
\label{sob.ineq}
\bigg ( \int_M |v|^{\frac{2n}{n-1}} \, \omega^n \bigg )^{\frac{n-1}{n}} \leq C \int_M |dv|_g^2 \, \omega^n + C \int_M |v|^2 \, \omega^n 
\end{equation}
for some uniform constant $C$.
\end{proposition}

\textbf{Proof.} 
We first consider the case when $v$ is supported in a coordinate neighborhood of some point in $\Sigma$. To that end, we fix a point on $\Sigma$, and let $(z_1,\hdots,z_{n-1},\zeta)$ be complex coordinates around that point. Let $g$ denote the Riemannian metric associated with the K\"ahler form $\omega$. Then the pull-back of the metric $g$ under the map $\Phi: (z_1,\hdots,z_{n-1},\xi) \mapsto (z_1,\hdots,z_{n-1},|\xi|^{\frac{1}{\beta}-1} \, \xi)$ is uniformly equivalent to the Euclidean metric. Therefore, the Sobolev inequality (\ref{sob.ineq}) holds if the function $v$ is suported in that coordinate chart. The general case follows in the standard way by using a partition of unity. \\

Finally, the Fredholm alternative established on page 20 of Donaldson's paper gives: 

\begin{theorem}[S.~Donaldson \cite{Donaldson}]
\label{fredholm}
Let $M$ be a compact K\"ahler manifold, and let $\Sigma$ be a complex hypersurface in $M$. Moreover, let $\omega$ be the background K\"ahler metric constructed above, and let $\hat{\omega}$ be a K\"ahler metric of class $\mathcal{C}^{,\alpha,\beta}$ which is uniformly equivalent to the background metric $\omega$. Then the operator 
\[\mathcal{C}^{2,\alpha,\beta} \to \mathcal{C}^{,\alpha,\beta}, \quad v \mapsto \Delta_{\hat{\omega}} v\] 
is Fredholm with Fredholm index zero.
\end{theorem}

\section{A $C^0$-estimate for solutions of $(\star_t)$}

\label{C0.estimate}

In this section, we establish a uniform $C^0$-estimate for solutions of $(\star_t)$. Throughout this section, we consider a pair $(u,c) \in \mathcal{C}^{2,\alpha,\beta} \times \mathbb{R}$ which satisfies the equation 
\begin{equation} 
\tag{$\star_t$} (\omega + \sqrt{-1} \, \partial\bar{\partial} u)^n = e^{tF-c} \, \omega^n
\end{equation}
for some $t \in [0,1]$. Since $u \in \mathcal{C}^{2,\alpha,\beta}$, it follows from standard elliptic regularity theory that $u \in C^\infty(M \setminus \Sigma)$. For abbreviation, let $\hat{\omega} = \omega + \sqrt{-1} \, \partial\bar{\partial} u$. 

Note that the function $u$ is only defined up to constants. We may normalize $u$ such that 
\begin{equation} 
\label{normalization}
\int_M u \, \omega_0^n = 0. 
\end{equation}
In the first step, we show that the function $u$ is uniformly bounded from above. To that end, we denote by $\Gamma_p: M \setminus \{p\} \to \mathbb{R}$ the Green's function associated with the operator $\Delta_{\omega_0}$ with pole at $p$. We may assume that the function $\Gamma_p$ is negative everywhere. 

\begin{proposition} 
\label{upper.bound}
We have 
\[u(p) \leq -n \int_M \Gamma_p \, \omega_0^n + \frac{\int_M \lambda \, |s|_h^{2\beta} \, \omega_0^n}{\int_M \omega_0^n}\] 
for all points $p \in M \setminus \Sigma$. In particular, we can find a uniform constant $N$ such that $N-u \geq 1$.
\end{proposition}

\textbf{Proof.} 
For abbreviation, let $v = u + \lambda \, |s|_h^{2\beta}$. Let us fix a point $p \in M \setminus \Sigma$, and let $\chi: M \to [0,1]$ be a smooth function which vanishes in a neighborhood of $\Sigma$. Using Green's formula, we obtain 
\begin{align*} 
\chi(p) \, v(p) 
&= \int_M \Gamma_p \, \Delta_{\omega_0}(\chi \, v) \, \omega_0^n + \frac{\int_M \chi \, v \, \omega_0^n}{\int_M \omega_0^n} \\ 
&= -\int_M \Gamma_p \, v \, \Delta_{\omega_0} \chi \, \omega_0^n - 2 \int_M v \, \langle d\chi,d\Gamma_p \rangle_{g_0} \, \omega_0^n \\ 
&+ \int_M \Gamma_p \, \chi \, \Delta_{\omega_0} v \, \omega_0^n + \frac{\int_M \chi \, v \, \omega_0^n}{\int_M \omega_0^n}. 
\end{align*} 
Using the identity $\omega_0 + \sqrt{-1} \, \partial\bar{\partial} v = \hat{\omega}$, we obtain $n + \Delta_{\omega_0} v \geq 0$. This implies 
\begin{align*} 
\chi(p) \, v(p) 
&\leq -\int_M \Gamma_p \, v \, \Delta_{\omega_0} \chi \, \omega_0^n - 2 \int_M v \, \langle d\chi,d\Gamma_p \rangle_{g_0} \, \omega_0^n \\ 
&- n \int_M \Gamma_p \, \chi \, \omega_0^n + \frac{\int_M \chi \, v \, \omega_0^n}{\int_M \omega_0^n}. 
\end{align*} 
We now choose the cutoff function $\chi$ such that $\chi=1$ outside a small neighborhood of $\Sigma$. We can choose $\chi$ in such a way that $\int_M (|\Delta_{\omega_0} \chi| + |d\chi|_{g_0}) \, \omega_0^n$ is arbitrarily small. Passing to the limit, we obtain 
\[v(p) \leq -n \int_M \Gamma_p \, \omega_0^n + \frac{\int_M v \, \omega_0^n}{\int_M \omega_0^n}.\] 
Since $u(p) \leq v(p)$ and $\int_M u \, \omega_0^n = 0$, the assertion follows. \\

\begin{lemma}
\label{L1.bound}
We have 
\[\int_M (N-u) \, \omega^n \leq C\] 
for some uniform constant $C$.
\end{lemma}

\textbf{Proof.} 
Since $u \in \mathcal{C}^{2,\alpha,\beta}$, we have $\sup_{M \setminus \Sigma} |du|_g < \infty$. This implies $\sup_{M \setminus \Sigma} |s|_h^{1-\beta} \, |du|_{g_0} < \infty$. Hence, integration by parts gives 
\begin{align*} 
&\int_{\{|s|_h \geq r\}} (N-u) \, \Delta_{\omega_0}(|s|_h^{2\beta}) \, \omega_0^n \\ 
&= -\int_{\{|s|_h \geq r\}} |s|_h^{2\beta} \, \Delta_{\omega_0} u \, \omega_0^n + O(r^{2\beta}). 
\end{align*} 
Using the inequality $n + \Delta_{\omega_0} (u + \lambda \, |s|_h^{2\beta}) \geq 0$, we obtain 
\begin{align*} 
&\int_{\{|s|_h \geq r\}} (N-u) \, \Delta_{\omega_0}(|s|_h^{2\beta}) \, \omega_0^n \\ 
&\leq \int_{\{|s|_h \geq r\}} |s|_h^{2\beta} \, \big ( n + \lambda \, \Delta_{\omega_0}(|s|_h^{2\beta}) \big ) \, \omega_0^n + O(r^{2\beta}). 
\end{align*} 
Taking the limit as $r \to 0$ gives 
\[\int_M (N-u) \, \Delta_{\omega_0}(|s|_h^{2\beta}) \, \omega_0^n \leq \int_M |s|_h^{2\beta} \, \big ( n + \lambda \, \Delta_{\omega_0}(|s|_h^{2\beta}) \big ) \, \omega_0^n \leq C.\] 
We next observe that 
\[\Delta_{\omega_0}(|s|_h^{2\beta}) \geq \delta \, |s|_h^{2\beta-2} - C_0\] 
for suitable constants $\delta,C_0 > 0$. Since $N-u$ is a positive function, we conclude that 
\[\int_M (N-u) \, (\delta \, |s|_h^{2\beta-2} - C_0) \, \omega_0^n \leq C.\] 
Using the normalization (\ref{normalization}), we obtain 
\[\int_M (N-u) \, |s|_h^{2\beta-2} \, \omega_0^n \leq C,\] 
hence 
\[\int_M (N-u) \, \omega^n \leq C.\] 
This completes the proof. \\

We next establish upper and lower bounds for the number $c$.

\begin{lemma}
\label{constant.c}
We have 
\[\int_M (e^{tF-c}-1) \, \omega^n = 0.\] 
In particular, $\inf_M tF \leq c \leq \sup_M tF$.
\end{lemma}

\textbf{Proof.} 
Let $\chi: M \to [0,1]$ be a smooth cutoff function which vanishes in a neighborhood of $\Sigma$. Then 
\begin{align*} 
\int_M \chi \, (e^{tF-c}-1) \, \omega^n 
&= \int_M \chi \, (\hat{\omega}^n - \omega^n) \\ 
&= \sum_{m=0}^{n-1} \int_M \chi \, (\hat{\omega} - \omega) \wedge \omega^m \wedge \hat{\omega}^{n-m-1} \\ 
&= \sum_{m=0}^{n-1} \int_M \chi \, \sqrt{-1} \, \partial\bar{\partial} u \wedge \omega^m \wedge \hat{\omega}^{n-m-1} \\ 
&= -\sum_{m=0}^{n-1} \int_M \sqrt{-1} \, \partial \chi \wedge \bar{\partial} u \wedge \omega^m \wedge \hat{\omega}^{n-m-1}. 
\end{align*}
Note that $\sup_{M \setminus \Sigma} |du|_g < \infty$ and $\sup_{M \setminus \Sigma} |\partial\bar{\partial} u|_g < \infty$ since $u \in \mathcal{C}^{2,\alpha,\beta}$. We now choose $\chi$ such that $\chi=1$ outside a small neighborhood of $\Sigma$. Moreover, we can arrange for $\int_M |d\chi|_g \, \omega^n$ to be arbitrarily small. Passing to the limit, the assertion follows. \\

In the next step, we establish an upper bound for the function $N-u$. This follows from an adaptation of the arguments in \cite{Tian}, Section 5.1. 

\begin{lemma} 
\label{iteration}
For each $p > 1$, we have 
\[\frac{4(p-1)}{p^2} \int_M \big | d \big ( (N-u)^{p/2} \big ) \big |_g^2 \, \omega^n \leq e^{\sup_M (tF) - \inf_M (tF)} \int_M (N-u)^p \, \omega^n.\] 
\end{lemma}

\textbf{Proof.} 
Let $\chi: M \to [0,1]$ be a smooth cutoff function which vanishes in a neighborhood of $\Sigma$. Then 
\begin{align*} 
&\int_M \chi \, (N-u)^{p-1} \, (e^{tF-c}-1) \, \omega^n \\ 
&= \int_M \chi \, (N-u)^{p-1} \, (\hat{\omega}^n - \omega^n) \\ 
&= \sum_{m=0}^{n-1} \int_M \chi \, (N-u)^{p-1} \, (\hat{\omega} - \omega) \wedge \omega^m \wedge \hat{\omega}^{n-m-1} \\ 
&= \sum_{m=0}^{n-1} \int_M \chi \, (N-u)^{p-1} \, \sqrt{-1} \, \partial\bar{\partial} u \wedge \omega^m \wedge \hat{\omega}^{n-m-1} \\ 
&= -\sum_{m=0}^{n-1} \int_M \sqrt{-1} \, \partial \big ( \chi \, (N-u)^{p-1} \big ) \wedge \bar{\partial} u \wedge \omega^m \wedge \hat{\omega}^{n-m-1} \\ 
&= (p-1) \sum_{m=0}^{n-1} \int_M \chi \, (N-u)^{p-2} \, \sqrt{-1} \, \partial u \wedge \bar{\partial} u \wedge \omega^m \wedge \hat{\omega}^{n-m-1} \\ 
&- \sum_{m=0}^{n-1} \int_M (N-u)^{p-1} \, \sqrt{-1} \, \partial \chi \wedge \bar{\partial} u \wedge \omega^m \wedge \hat{\omega}^{n-m-1}. 
\end{align*} 
Note that $\sup_{M \setminus \Sigma} |du|_g < \infty$ and $\sup_{M \setminus \Sigma} |\partial\bar{\partial} u|_g < \infty$ since $u \in \mathcal{C}^{2,\alpha,\beta}$. As above, we may choose $\chi$ such that $\chi=1$ outside a small neighborhood of $\Sigma$ and $\int_M |d\chi|_g \, \omega^n$ is arbitrarily small. Passing to the limit, we obtain 
\begin{align*} 
&\int_M (N-u)^{p-1} \, (e^{tF-c}-1) \, \omega^n \\ 
&= (p-1) \sum_{m=0}^{n-1} \int_M (N-u)^{p-2} \, \sqrt{-1} \, \partial u \wedge \bar{\partial} u \wedge \omega^m \wedge \hat{\omega}^{n-m-1}. 
\end{align*} 
This implies 
\begin{align*} 
&\int_M (N-u)^{p-1} \, (e^{tF-c}-1) \, \omega^n \\ 
&\geq (p-1) \int_M (N-u)^{p-2} \, \sqrt{-1} \, \partial u \wedge \bar{\partial} u \wedge \omega^{n-1} \\ 
&= (p-1) \int_M (N-u)^{p-2} \, |du|_g^2 \, \omega^n \\ 
&= \frac{4(p-1)}{p^2} \int_M \big | d \big ( (N-u)^{p/2} \big ) \big |_g^2 \, \omega^n. 
\end{align*} 
On the other hand, we have $c \geq \inf_M (tF)$ by Lemma \ref{constant.c}. From this, we deduce that 
\begin{align*} 
\int_M (N-u)^{p-1} \, (e^{tF-c}-1) \, \omega^n 
&\leq e^{\sup_M (tF) - \inf_M (tF)} \int_M (N-u)^{p-1} \, \omega^n \\ 
&\leq e^{\sup_M (tF) - \inf_M (tF)} \int_M (N-u)^p \, \omega^n. 
\end{align*}
Putting these facts together, the assertion follows. \\

\begin{lemma}
\label{Lp.bound}
We have 
\[\int_M (N-u)^{\frac{2n}{n-1}} \, \omega^n \leq C\] 
for some uniform constant $C$. 
\end{lemma}

\textbf{Proof.} 
Applying Lemma \ref{iteration} with $p=2$ gives 
\[\int_M |d(N-u)|_g^2 \, \omega^n \leq e^{\sup_M (tF) - \inf_M (tF)} \int_M (N-u)^2 \, \omega^n.\] 
Moreover, we have 
\[\bigg ( \int_M (N-u)^{\frac{2n}{n-1}} \, \omega^n \bigg )^{\frac{n-1}{n}} \leq C \int_M |d(N-u)|_g^2 \, \omega^n + C \int_M (N-u)^2 \, \omega^n\] 
by Proposition \ref{sobolev}. Putting these facts together, we obtain 
\begin{align*} 
&\bigg ( \int_M (N-u)^{\frac{2n}{n-1}} \, \omega^n \bigg )^{\frac{n-1}{n}} \\ 
&\leq C \int_M (N-u)^2 \, \omega^n \leq C \, \bigg ( \int_M (N-u) \, \omega^n \bigg )^{\frac{2}{n+1}} \, \bigg ( \int_M (N-u)^{\frac{2n}{n-1}} \, \omega^n \bigg )^{\frac{n-1}{n+1}} 
\end{align*} 
for some uniform constant $C$. Thus, we conclude that 
\[\bigg ( \int_M (N-u)^{\frac{2n}{n-1}} \, \omega^n \bigg )^{\frac{n-1}{2n}} \leq C \int_M (N-u) \, \omega^n.\] 
Hence, the assertion follows from Lemma \ref{L1.bound}. \\

\begin{proposition} 
\label{lower.bound}
We have $\sup_M (N-u) \leq C$ for some uniform constant $C$.
\end{proposition}

\textbf{Proof.} 
By Lemma \ref{Lp.bound}, we have 
\[\int_M (N-u)^{\frac{2n}{n-1}} \, \omega^n \leq C.\] 
Moreover, using Proposition \ref{sobolev} and Lemma \ref{iteration}, we obtain 
\begin{align*} 
&\bigg ( \int_M (N-u)^{\frac{pn}{n-1}} \, \omega^n \bigg )^{\frac{n-1}{n}} \\ 
&\leq C \int_M \big | d((N-u)^{p/2}) \big |^2 \, \omega^n + C \int_M (N-u)^p \, \omega^n \\ 
&\leq C \, \Big ( \frac{p^2}{4(p-1)} \, e^{\sup_M (tF) - \inf_M (tF)} + 1 \Big ) \int_M (N-u)^p \, \omega^n
\end{align*}
for each $p > 1$. Hence, it follows from the standard Moser iteration technique that $\sup_M (N-u) \leq C$ for some uniform constant $C$. \\

Combining Proposition \ref{upper.bound} and Proposition \ref{lower.bound}, we conclude that $\sup_M |u| \leq C$ for some uniform constant $C$.

\section{An estimate for $\partial \bar{\partial} u$}

\label{C2.estimate}

In this section, we establish a uniform estimate for the Laplacian $\Delta_\omega u$. The following lemma will be useful:

\begin{lemma}
\label{barrier}
If $\varepsilon > 0$ is sufficiently small, then 
\[\sqrt{-1} \, \partial\bar{\partial}(|s|_h^{2\varepsilon}) + \omega \geq \varepsilon^2 \, |s|_h^{2\varepsilon} \, \sqrt{-1} \, \partial \log(|s|_h^2) \wedge \bar{\partial} \log(|s|_h^2) \geq 0.\] 
\end{lemma}

\textbf{Proof.} 
We can find a positive constant $C$ such that 
\[\sqrt{-1} \, \partial\bar{\partial} \log(|s|_h^2) + C \, \omega_0 \geq 0.\] 
This implies 
\[\sqrt{-1} \, \partial\bar{\partial}(|s|_h^{2\varepsilon}) + C \, \varepsilon \, |s|_h^{2\varepsilon} \, \omega_0 \geq \varepsilon^2 \, |s|_h^{2\varepsilon} \, \sqrt{-1} \, \partial \log(|s|_h^2) \wedge \bar{\partial} \log(|s|_h^2).\] 
If we choose $\varepsilon$ sufficiently small, then $C \, \varepsilon \, \omega_0 \leq \omega$. From this, the assertion follows. \\

The proof of the Laplacian estimate largely follows the arguments in Yau's paper \cite{Yau2}. There is an additional technical issue in that function $e^{-Lu} \, (n + \Delta_\omega u)$ may attain its maximum on the singular set $\Sigma$. This obstacle can be overcome using a trick due to T.~Jeffres (cf. \cite{Jeffres2}, Section 4). Jeffres studied the case of K\"ahler-Einstein metrics with negative scalar curvature, and proved a $C^0$ estimate in that setting. A similar trick works for the Laplacian estimate.

\begin{proposition}
\label{c2}
Let $u \in \mathcal{C}^{2,\alpha,\beta}$ be a solution of $(\star_t)$ for some $t \in [0,1]$. Then $\Delta_\omega u \leq C$ for some uniform constant $C$.
\end{proposition}

\textbf{Proof.} 
Let us fix a point $q \in M \setminus \Sigma$. We may choose complex coordinates around $q$ so that $\omega_{\mu\bar{\nu}} = \delta_{\mu\nu}$ and $u_{\mu\bar{\nu}} = u_{\mu\bar{\mu}} \, \delta_{\mu\nu}$. By Proposition \ref{curvature}, the curvature tensor of the background metric $\omega$ is uniformly bounded. Hence, we can find a positive constant $L$ such that 
\[L + \inf_{\mu\neq\nu} R_{\mu\bar{\mu}\nu\bar{\nu}} \geq 2.\] 
Using equation (2.18) in \cite{Yau2}, we obtain 
\begin{align*} 
&\Delta_{\hat{\omega}}(e^{-Lu} \, (n + \Delta_\omega u)) \\ 
&\geq e^{-Lu} \, \Big ( \Delta_\omega(tF) - n^2 \inf_{\mu\neq\nu} R_{\mu\bar{\mu}\nu\bar{\nu}} \Big ) - nL \, e^{-Lu} \, (n + \Delta_\omega u) \\ 
&+ \Big ( L + \inf_{\mu\neq\nu} R_{\mu\bar{\mu}\nu\bar{\nu}} \Big ) \, e^{-Lu} \, (n + \Delta_\omega u) \, \sum_{\mu=1}^n \frac{1}{1+u_{\mu\bar{\mu}}} 
\end{align*}
at the point $q$. By Proposition \ref{properties.F}, we have $|\Delta_\omega F| \leq C$. Moreover, Proposition \ref{curvature} implies that the background metric $\omega$ has bounded curvature. This implies 
\begin{align*} 
\Delta_{\hat{\omega}}(e^{-Lu} \, (n + \Delta_\omega u)) 
&\geq 2 \, e^{-Lu} \, (n + \Delta_\omega u) \, \sum_{\mu=1}^n \frac{1}{1+u_{\mu\bar{\mu}}} \\ 
&-nL \, e^{-Lu} \, (n + \Delta_\omega u) - e^{-Lu} \, C 
\end{align*} 
at the point $q$. We next observe that 
\[n + \Delta_\omega u \geq n \, \bigg ( \prod_{\mu=1}^n (1+u_{\mu\bar{\mu}}) \bigg )^{\frac{1}{n}} = n \, e^{\frac{tF-c}{n}}\] 
and 
\[\sum_{\mu=1}^n \frac{1}{1+u_{\mu\bar{\mu}}} \geq \bigg ( \frac{n + \Delta_\omega u}{\prod_{\mu=1}^n (1+u_{\mu\bar{\mu}})} \bigg )^{\frac{1}{n-1}} = e^{-\frac{tF-c}{n-1}} \, (n + \Delta_\omega u)^{\frac{1}{n-1}}\] 
(cf. \cite{Yau2}, equation (2.19)). Thus, we conclude that 
\begin{align*} 
&\Delta_{\hat{\omega}}(e^{-Lu} \, (n + \Delta_\omega u)) \\ 
&\geq n \, e^{\frac{tF-c}{n}} \, e^{-Lu} \, \sum_{\mu=1}^n \frac{1}{1+u_{\mu\bar{\mu}}} + e^{-\frac{tF-c}{n-1}} \, e^{-Lu} \, (n + \Delta_\omega u)^{\frac{n}{n-1}} \\ 
&-nL \, e^{-Lu} \, (n + \Delta_\omega u) - e^{-Lu} \, C. 
\end{align*} 
Moreover, it follows from Propositions \ref{upper.bound} and \ref{lower.bound} that $|u| \leq C$ for some uniform constant $C$. Putting these facts together, we obtain 
\begin{align*} 
&\Delta_{\hat{\omega}}(e^{-Lu} \, (n + \Delta_\omega u)) \\ 
&\geq \kappa \sum_{\mu=1}^n \frac{1}{1+u_{\mu\bar{\mu}}} + \kappa \, e^{-\frac{nLu}{n-1}} \, (n + \Delta_\omega u)^{\frac{n}{n-1}} - nL \, e^{-Lu} \, (n + \Delta_\omega u) - C 
\end{align*} 
for uniform constants $\kappa, C > 0$.

Let us fix $\varepsilon > 0$ sufficiently small. By Lemma \ref{barrier}, we have 
\[\sqrt{-1} \, \partial\bar{\partial} (|s|_h^{2\varepsilon}) + \omega \geq 0,\] 
hence 
\[\Delta_{\hat{\omega}} (|s|_h^{2\varepsilon}) + \sum_{\mu=1}^n \frac{1}{1+u_{\mu\bar{\mu}}} \geq 0.\] 
Consequently, the function 
\[H = e^{-Lu} \, (n + \Delta_\omega u) + \kappa \, |s|_h^{2\varepsilon}\] 
satisfies 
\begin{align*} 
\Delta_{\hat{\omega}} H 
&\geq \kappa \, \Delta_{\hat{\omega}}(|s|_h^{2\varepsilon}) + \kappa \sum_{\mu=1}^n \frac{1}{1+u_{\mu\bar{\mu}}} + \kappa \, \big ( H - \kappa \, |s|_h^{2\varepsilon} \big )^{\frac{n}{n-1}} \\ 
&- nL \, \big ( H - \kappa \, |s|_h^{2\varepsilon} \big ) - C \\ 
&\geq \kappa \, \big ( H - \kappa \, |s|_h^{2\varepsilon} \big )^{\frac{n}{n-1}} - nL \, \big ( H - \kappa \, |s|_h^{2\varepsilon} \big ) - C. 
\end{align*} 
Since the function $u$ is of class $\mathcal{C}^{2,\alpha,\beta}$, the function $\Delta_\omega u$ is of class $\mathcal{C}^{,\alpha,\beta}$. Consequently, if $\varepsilon > 0$ is small, then the function $H$ attains its maximum at some point $q \in M \setminus \Sigma$. The maximum principle then implies that $\sup_{M \setminus \Sigma} H \leq C$ for some constant $C$. Consequently, the function $\Delta_\omega u$ is uniformly bounded from above, as claimed. \\

It follows from $(\star_t)$ that the function $\frac{\hat{\omega}^n}{\omega^n}$ is uniformly bounded from above and below. Using Proposition \ref{c2}, we conclude that $a_1 \, \omega \leq \hat{\omega} \leq a_2 \, \omega$ for uniform constants $a_1,a_2 > 0$.

\section{An estimate for the covariant derivative of $\partial\bar{\partial} u$}

\label{C3.estimate} 
 
As above, we consider a pair $(u,c) \in \mathcal{C}^{2,\alpha,\beta} \times \mathbb{R}$ which satisfies the equation 
\begin{equation} 
\tag{$\star_t$} (\omega + \sqrt{-1} \, \partial\bar{\partial} u)^n = e^{tF-c} \, \omega^n
\end{equation}
for some $t \in [0,1]$. Moreover, we let $\hat{\omega} = \omega + \sqrt{-1} \, \partial\bar{\partial} u$. 

Our goal in this section is to establish a uniform bound for the covariant derivative of $\partial\bar{\partial} u$. This estimate will be proved using the maximum principle. However, in order to apply the maximum principle, it is necessary to analyze the behavior of the third derivatives of $u$ near $\Sigma$. To that end, we fix a point $p \in \Sigma$. Moreover, let $(z_1,\hdots,z_{n-1},\zeta)$ be complex coordinates around $p$.

\begin{lemma}
\label{third.order.derivatives.a}
We have 
\begin{align*}
&\frac{\partial^3 u}{\partial z_i \, \partial \bar{z}_j \, \partial z_k} \in \mathcal{C}^{,\alpha,\beta}, \\ 
&|\zeta|^{1-\beta} \, \frac{\partial^3 u}{\partial z_i \, \partial \bar{\zeta} \, \partial z_k} \in \mathcal{C}_0^{,\alpha,\beta}, \\ 
&|\zeta|^{1-\beta} \, \frac{\partial^3 u}{\partial \zeta \, \partial \bar{z}_j \, \partial z_k} \in \mathcal{C}_0^{,\alpha,\beta}, \\
&|\zeta|^{2-2\beta} \, \frac{\partial^3 u}{\partial \zeta \, \partial \bar{\zeta} \, \partial z_k} \in \mathcal{C}^{,\alpha,\beta} 
\end{align*} 
for all $i,j,k \in \{1,\hdots,n-1\}$.
\end{lemma}

\textbf{Proof.} 
Let us fix an integer $k \in \{1,\hdots,n-1\}$. Using the formula 
\[\omega = \omega_0 + \lambda \, \sqrt{-1} \, \partial\bar{\partial}(\rho \, |\zeta|^{2\beta}),\] 
we obtain 
\[\frac{\partial}{\partial z_k} \omega = \frac{\partial}{\partial z_k} \omega_0 + \lambda \, \sqrt{-1} \, \partial\bar{\partial} \Big ( \frac{\partial \rho}{\partial z_k} \, |\zeta|^{2\beta} \Big ).\] 
Hence, if $\alpha>0$ is sufficiently small, then the $(1,1)$-form $\frac{\partial}{\partial z_k} \omega$ is of class $\mathcal{C}^{,\alpha,\beta}$. This implies $\text{\rm tr}_\omega \big ( \frac{\partial}{\partial z_k} \omega \big ) \in \mathcal{C}^{,\alpha,\beta}$. Since $\hat{\omega}$ is uniformly equivalent to $\omega$ and of class $\mathcal{C}^{,\alpha,\beta}$, we also have $\text{\rm tr}_{\hat{\omega}} \big ( \frac{\partial}{\partial z_k} \omega \big ) \in \mathcal{C}^{,\alpha,\beta}$. Moreover, it follows from (\ref{formula.F}) that $\frac{\partial}{\partial z_k} F \in \mathcal{C}^{,\alpha,\beta}$. We now define 
\[v = \frac{\partial}{\partial z_k} u \in \mathcal{C}^{,\alpha,\beta}\] 
and 
\[f = \text{\rm tr}_\omega \Big ( \frac{\partial}{\partial z_k} \omega \Big ) - \text{\rm tr}_{\hat{\omega}} \Big ( \frac{\partial}{\partial z_k} \omega \Big ) + \frac{\partial}{\partial z_k} (tF) \in \mathcal{C}^{,\alpha,\beta}.\] 
Differentiating the equation $(\star_t)$ with respect to $z_k$, we conclude that $\Delta_{\hat{\omega}} v = f$ away from $\Sigma$. Hence, Corollary \ref{regularity.3} implies that $v \in \mathcal{C}^{2,\alpha,\beta}$. From this, the assertion follows. \\

\begin{lemma}
\label{third.order.derivatives.b}
We have 
\[|\zeta|^{2-2\beta} \, \Big ( \frac{\partial^3 u}{\partial \zeta \, \partial \bar{z}_j \, \partial \zeta} + \frac{1-\beta}{\zeta} \, \frac{\partial^2 u}{\partial \zeta \, \partial \bar{z}_j} \Big ) = O(|\zeta|^{\alpha\beta})\] 
for all $j \in \{1,\hdots,n-1\}$.
\end{lemma}

\textbf{Proof.} 
Let us fix an integer $j \in \{1,\hdots,n-1\}$ and let $v = \frac{\partial}{\partial \bar{z}_j} u$. By Lemma \ref{third.order.derivatives.a}, we have 
\[|\zeta|^{2-2\beta} \, \frac{\partial^2 v}{\partial \zeta \, \partial \bar{\zeta}} \in \mathcal{C}^{,\alpha,\beta}.\] 
Applying Proposition \ref{auxiliary.result} to the function $\zeta \mapsto v(z_1,\hdots,z_{n-1},\zeta)$, we conclude that 
\[|\zeta|^{2-2\beta} \, \Big ( \frac{\partial^2 v}{\partial \zeta \, \partial \zeta} + \frac{1-\beta}{\zeta} \, \frac{\partial v}{\partial \zeta} \Big ) = O(|\zeta|^{\alpha\beta}),\] 
as claimed. \\

In the next step, we study the third order covariant derivatives of $u$ with respect to the background edge metric $\omega$. We first analyze the Christoffel symbols of the background metric $g$.

\begin{lemma}
\label{christoffel}
The Christoffel symbols of $g$ satisfy 
\[\Gamma_{ik}^j \in \mathcal{C}^{,\alpha,\beta}, \quad |\zeta|^{1-\beta} \, \Gamma_{i\zeta}^k \in \mathcal{C}_0^{,\alpha,\beta}, \quad |\zeta|^{\beta-1} \, \Gamma_{ik}^\zeta \in \mathcal{C}_0^{,\alpha,\beta}\] 
and 
\[\Gamma_{i\zeta}^\zeta \in \mathcal{C}^{,\alpha,\beta}, \quad |\zeta|^{2-2\beta} \, \Gamma_{\zeta\zeta}^k \in \mathcal{C}_0^{,\alpha,\beta}, \quad |\zeta|^{1-\beta} \, \Big ( \Gamma_{\zeta\zeta}^\zeta + \frac{1-\beta}{\zeta} \Big ) \in \mathcal{C}_0^{,\alpha,\beta}\] 
for all $i,j,k \in \{1,\hdots,n-1\}$.
\end{lemma}

\textbf{Proof.} 
Recall that 
\begin{align*} 
\omega 
&= \omega_0 + \lambda \, |\zeta|^{2\beta} \, \sqrt{-1} \, \partial\bar{\partial} \rho \\ 
&+ \lambda\beta \, |\zeta|^{2\beta-2} \, \zeta \, \sqrt{-1} \, \partial \rho \wedge d\bar{\zeta} \\ 
&+ \lambda\beta \, |\zeta|^{2\beta-2} \, \bar{\zeta} \, \sqrt{-1} \, d\zeta \wedge \bar{\partial} \rho \\ 
&+ \lambda\beta^2\rho \, \sqrt{-1} \, |\zeta|^{2\beta-2} \, d\zeta \wedge d\bar{\zeta}. 
\end{align*}
Hence, if $\alpha>0$ is sufficiently small, then we have 
\[g^{k\bar{l}} \in \mathcal{C}^{,\alpha,\beta}, \quad |\zeta|^{\beta-1} \, g^{k\bar{\zeta}} \in \mathcal{C}_0^{,\alpha,\beta}, \quad |\zeta|^{2\beta-2} \, g^{\zeta\bar{\zeta}} \in \mathcal{C}^{,\alpha,\beta}.\] 
Moreover, since $\beta \in (0,\frac{1}{2})$, we have 
\begin{align*} 
&\frac{\partial}{\partial z_i} g \Big ( \frac{\partial}{\partial z_k},\frac{\partial}{\partial \bar{z}_l} \Big ) \in \mathcal{C}^{,\alpha,\beta}, \\ 
&|\zeta|^{1-\beta} \, \frac{\partial}{\partial z_i} g \Big ( \frac{\partial}{\partial \zeta},\frac{\partial}{\partial \bar{z}_l} \Big ) \in \mathcal{C}_0^{,\alpha,\beta}, \\ 
&|\zeta|^{2-2\beta} \, \frac{\partial}{\partial z_i} g \Big ( \frac{\partial}{\partial \zeta},\frac{\partial}{\partial \bar{\zeta}} \Big ) \in \mathcal{C}^{,\alpha,\beta}, \\ 
&|\zeta|^{2-2\beta} \, \Big [ \frac{\partial}{\partial \zeta} g \Big ( \frac{\partial}{\partial \zeta},\frac{\partial}{\partial \bar{z}_l} \Big ) + \frac{1-\beta}{\zeta} \, g \Big ( \frac{\partial}{\partial \zeta},\frac{\partial}{\partial \bar{z}_l} \Big ) \Big ] \in \mathcal{C}_0^{,\alpha,\beta}, \\ 
&|\zeta|^{3-3\beta} \, \Big [ \frac{\partial}{\partial \zeta} g \Big ( \frac{\partial}{\partial \zeta},\frac{\partial}{\partial \bar{\zeta}} \Big ) + \frac{1-\beta}{\zeta} \, g \Big ( \frac{\partial}{\partial \zeta},\frac{\partial}{\partial \bar{\zeta}} \Big ) \Big ] \in \mathcal{C}_0^{,\alpha,\beta}
\end{align*}
if $\alpha>0$ is sufficiently small. This implies 
\begin{align*}
&\Gamma_{ik}^j = \sum_{l=1}^{n-1} g^{j\bar{l}} \, \partial_i g_{k\bar{l}} + g^{j\bar{\zeta}} \, \partial_i g_{k\bar{\zeta}} \in \mathcal{C}^{,\alpha,\beta}, \\ 
&|\zeta|^{1-\beta} \, \Gamma_{i\zeta}^k = |\zeta|^{1-\beta} \sum_{l=1}^{n-1} g^{k\bar{l}} \, \partial_i g_{\zeta\bar{l}} + |\zeta|^{1-\beta} \, g^{k\bar{\zeta}} \, \partial_i g_{\zeta\bar{\zeta}} \in \mathcal{C}_0^{,\alpha,\beta}, \\ 
&|\zeta|^{\beta-1} \, \Gamma_{ik}^\zeta = |\zeta|^{\beta-1}  \sum_{l=1}^{n-1} g^{\zeta\bar{l}} \, \partial_i g_{k\bar{l}} + |\zeta|^{\beta-1} \, g^{\zeta\bar{\zeta}} \, \partial_i g_{k\bar{\zeta}} \in \mathcal{C}_0^{,\alpha,\beta}, \\ 
&\Gamma_{i\zeta}^\zeta = \sum_{l=1}^{n-1} g^{\zeta\bar{l}} \, \partial_i g_{\zeta\bar{l}} + g^{\zeta\bar{\zeta}} \, \partial_i g_{\zeta\bar{\zeta}} \in \mathcal{C}^{,\alpha,\beta}. 
\end{align*} 
Moreover, we have 
\begin{align*} 
|\zeta|^{2-2\beta} \, \Gamma_{\zeta\zeta}^k 
&= |\zeta|^{2-2\beta} \sum_{l=1}^{n-1} g^{k\bar{l}} \, \Big ( \partial_\zeta g_{\zeta\bar{l}} + \frac{1-\beta}{\zeta} \, g_{\zeta\bar{l}} \Big ) \\ 
&+ |\zeta|^{2-2\beta} \, g^{k\bar{\zeta}} \, \Big ( \partial_\zeta g_{\zeta\bar{\zeta}} + \frac{1-\beta}{\zeta} \, g_{\zeta\bar{\zeta}} \Big ) \in \mathcal{C}_0^{,\alpha,\beta} 
\end{align*} 
and 
\begin{align*} 
|\zeta|^{1-\beta} \, \Big ( \Gamma_{\zeta\zeta}^\zeta + \frac{1-\beta}{\zeta} \Big ) 
&= |\zeta|^{1-\beta} \sum_{l=1}^{n-1} g^{\zeta\bar{l}} \, \Big ( \partial_\zeta g_{\zeta\bar{l}} + \frac{1-\beta}{\zeta} \, g_{\zeta\bar{l}} \Big ) \\ 
&+ |\zeta|^{1-\beta} \, g^{\zeta\bar{\zeta}} \, \Big ( \partial_\zeta g_{\zeta\bar{\zeta}} + \frac{1-\beta}{\zeta} \, g_{\zeta\bar{\zeta}} \Big ) \in \mathcal{C}_0^{,\alpha,\beta}. 
\end{align*} 
This completes the proof. \\

\begin{proposition}
\label{covariant.derivatives.a}
The third order covariant derivatives of $u$ satisfy 
\begin{align*}
&(D^3 u) \Big ( \frac{\partial}{\partial z_i},\frac{\partial}{\partial \bar{z}_j}, \frac{\partial}{\partial z_k} \Big ) \in \mathcal{C}^{,\alpha,\beta}, \\ 
&|\zeta|^{1-\beta} \, (D^3 u) \Big ( \frac{\partial}{\partial z_i},\frac{\partial}{\partial \bar{\zeta}},\frac{\partial}{\partial z_k} \Big ) \in \mathcal{C}_0^{,\alpha,\beta}, \\ 
&|\zeta|^{1-\beta} \, (D^3 u) \Big ( \frac{\partial}{\partial \zeta},\frac{\partial}{\partial \bar{z}_j}, \frac{\partial}{\partial z_k} \Big ) \in \mathcal{C}_0^{,\alpha,\beta}, \\
&|\zeta|^{2-2\beta} \, (D^3 u) \Big ( \frac{\partial}{\partial \zeta},\frac{\partial}{\partial \bar{\zeta}},\frac{\partial}{\partial z_k} \Big ) \in \mathcal{C}^{,\alpha,\beta} 
\end{align*} 
for all $i,j,k \in \{1,\hdots,n-1\}$. Moreover, we have 
\[|\zeta|^{2-2\beta} \, (D^3 u) \Big ( \frac{\partial}{\partial \zeta}, \frac{\partial}{\partial \bar{z}_j},\frac{\partial}{\partial \zeta} \Big ) = O(|\zeta|^{\alpha\beta})\] 
for all $j \in \{1,\hdots,n-1\}$.
\end{proposition}

\textbf{Proof.} 
Using Lemma \ref{third.order.derivatives.a} and Lemma \ref{christoffel}, we obtain 
\begin{align*} 
&(D^3 u) \Big ( \frac{\partial}{\partial z_i},\frac{\partial}{\partial \bar{z}_j}, \frac{\partial}{\partial z_k} \Big ) \\ 
&= \frac{\partial^3 u}{\partial z_i \, \partial \bar{z}_j \, \partial z_k} - \sum_{l=1}^{n-1} \Gamma_{ik}^l \, \frac{\partial^2 u}{\partial z_l \, \partial \bar{z}_j} - \Gamma_{ik}^\zeta \, \frac{\partial^2 u}{\partial \zeta \, \partial \bar{z}_j} \in \mathcal{C}^{,\alpha,\beta}, \\[1.5mm] 
&|\zeta|^{1-\beta} \, (D^3 u) \Big ( \frac{\partial}{\partial z_i},\frac{\partial}{\partial \bar{\zeta}},\frac{\partial}{\partial z_k} \Big ) \\ 
&= |\zeta|^{1-\beta} \, \frac{\partial^3 u}{\partial z_i \, \partial \bar{\zeta} \, \partial z_k} - |\zeta|^{1-\beta} \sum_{l=1}^{n-1} \Gamma_{ik}^l \, \frac{\partial^2 u}{\partial z_l \, \partial \bar{\zeta}} - |\zeta|^{1-\beta} \, \Gamma_{ik}^\zeta \, \frac{\partial^2 u}{\partial \zeta \, \partial \bar{\zeta}} \in \mathcal{C}_0^{,\alpha,\beta}, \\[1.5mm] 
&|\zeta|^{1-\beta} \, (D^3 u) \Big ( \frac{\partial}{\partial \zeta},\frac{\partial}{\partial \bar{z}_j}, \frac{\partial}{\partial z_k} \Big ) \\ 
&= |\zeta|^{1-\beta} \, \frac{\partial^3 u}{\partial \zeta \, \partial \bar{z}_j \, \partial z_k} - |\zeta|^{1-\beta} \sum_{l=1}^{n-1} \Gamma_{k\zeta}^l \, \frac{\partial^2 u}{\partial z_l \, \partial \bar{z}_j} - |\zeta|^{1-\beta} \, \Gamma_{k\zeta}^\zeta \, \frac{\partial^2 u}{\partial \zeta \, \partial \bar{z}_j} \in \mathcal{C}_0^{,\alpha,\beta}, \\[1.5mm] 
&|\zeta|^{2-2\beta} \, (D^3 u) \Big ( \frac{\partial}{\partial \zeta},\frac{\partial}{\partial \bar{\zeta}},\frac{\partial}{\partial z_k} \Big ) \\ 
&= |\zeta|^{2-2\beta} \, \frac{\partial^3 u}{\partial \zeta \, \partial \bar{\zeta} \, \partial z_k} - |\zeta|^{2-2\beta} \sum_{l=1}^{n-1} \Gamma_{k\zeta}^l \, \frac{\partial^2 u}{\partial z_l \, \partial \bar{\zeta}} - |\zeta|^{2-2\beta} \, \Gamma_{k\zeta}^\zeta \, \frac{\partial^2 u}{\partial \zeta \, \partial \bar{\zeta}} \in \mathcal{C}^{,\alpha,\beta}. 
\end{align*} 
We next observe that 
\begin{align*} 
|\zeta|^{2-2\beta} \, (D^3 u) \Big ( \frac{\partial}{\partial \zeta}, \frac{\partial}{\partial \bar{z}_j},\frac{\partial}{\partial \zeta} \Big ) 
&= |\zeta|^{2-2\beta} \, \Big ( \frac{\partial^3 u}{\partial \zeta \, \partial \bar{z}_j \, \partial \zeta} + \frac{1-\beta}{\zeta} \, \frac{\partial^2 u}{\partial \zeta \, \partial \bar{z}_j} \Big ) \\ 
&- |\zeta|^{2-2\beta} \, \sum_{j=1}^{n-1} \Gamma_{\zeta\zeta}^l \, \frac{\partial^2 u}{\partial z_l \, \partial \bar{z}_j} \\ 
&- |\zeta|^{2-2\beta} \, \Big ( \Gamma_{\zeta\zeta}^\zeta + \frac{1-\beta}{\zeta} \Big ) \, \frac{\partial^2 u}{\partial \zeta \, \partial \bar{z}_j}. 
\end{align*}
Since 
\[|\zeta|^{2-2\beta} \, \Big ( \frac{\partial^3 u}{\partial \zeta \, \partial \bar{z}_j \, \partial \zeta} + \frac{1-\beta}{\zeta} \, \frac{\partial^2 u}{\partial \zeta \, \partial \bar{z}_j} \Big ) = O(|\zeta|^{\alpha\beta})\] 
and 
\[|\zeta|^{1-\beta} \, \Big ( \Gamma_{\zeta\zeta}^\zeta + \frac{1-\beta}{\zeta} \Big ) = O(|\zeta|^{\alpha\beta}),\] 
we conclude that 
\[|\zeta|^{2-2\beta} \, (D^3 u) \Big ( \frac{\partial}{\partial \zeta}, \frac{\partial}{\partial \bar{z}_j},\frac{\partial}{\partial \zeta} \Big ) = O(|\zeta|^{\alpha\beta}),\] 
as claimed. \\

\begin{proposition}
\label{covariant.derivatives.b}
We have 
\[|\zeta|^{3-3\beta} \, (D^3 u) \Big ( \frac{\partial}{\partial \zeta},\frac{\partial}{\partial \bar{\zeta}},\frac{\partial}{\partial \zeta} \Big ) = O(|\zeta|^{\alpha\beta}).\] 
\end{proposition}

\textbf{Proof.} 
Differentiating the Monge-Amp\`ere equation $(\star_t)$ gives 
\begin{equation}
\label{identity.1}
\sum_{\mu,\nu=1}^n \hat{g}^{\mu\bar{\nu}} \, (D^3 u) \Big ( \frac{\partial}{\partial \zeta},\frac{\partial}{\partial \bar{z}_\nu},\frac{\partial}{\partial z_\mu} \Big ) = \frac{\partial}{\partial \zeta}(tF). 
\end{equation}
Since the metric $\hat{\omega}$ is of class $\mathcal{C}^{,\alpha,\beta}$, we have 
\[\hat{g}^{k\bar{l}} \in \mathcal{C}^{,\alpha,\beta}, \quad |\zeta|^{\beta-1} \, \hat{g}^{k\bar{\zeta}} \in \mathcal{C}_0^{,\alpha,\beta}, \quad |\zeta|^{2\beta-2} \, \hat{g}^{\zeta\bar{\zeta}} \in \mathcal{C}^{,\alpha,\beta}.\] 
Using Proposition \ref{covariant.derivatives.a}, we obtain 
\begin{align*}
&|\zeta|^{1-\beta} \, \hat{g}^{k\bar{l}} \, (D^3 u) \Big ( \frac{\partial}{\partial \zeta},\frac{\partial}{\partial \bar{z}_l},\frac{\partial}{\partial z_k} \Big ) = O(|\zeta|^{\alpha\beta}), \\ 
&|\zeta|^{1-\beta} \, \hat{g}^{\zeta\bar{l}} \, (D^3 u) \Big ( \frac{\partial}{\partial \zeta},\frac{\partial}{\partial \bar{z}_l},\frac{\partial}{\partial \zeta} \Big ) = O(|\zeta|^{\alpha\beta}), \\ 
&|\zeta|^{1-\beta} \, \hat{g}^{k\bar{\zeta}} \, (D^3 u) \Big ( \frac{\partial}{\partial \zeta},\frac{\partial}{\partial \bar{\zeta}},\frac{\partial}{\partial z_k} \Big ) = O(|\zeta|^{\alpha\beta}) 
\end{align*} 
for $k,l \in \{1,\hdots,n-1\}$. Substituting this into (\ref{identity.1}), we obtain 
\begin{equation} 
\label{identity.2} 
|\zeta|^{1-\beta} \, \hat{g}^{\zeta\bar{\zeta}} \, (D^3 u) \Big (\frac{\partial}{\partial \zeta},\frac{\partial}{\partial \bar{\zeta}},\frac{\partial}{\partial \zeta} \Big ) = |\zeta|^{1-\beta} \, \frac{\partial}{\partial \zeta}(tF) + O(|\zeta|^{\alpha\beta}). 
\end{equation}
It follows from (\ref{formula.F}) that $|\zeta|^{1-\beta} \, \frac{\partial}{\partial \zeta} F = O(|\zeta|^\beta)$. Putting these facts together, we conclude that 
\[|\zeta|^{1-\beta} \, \hat{g}^{\zeta\bar{\zeta}} \, (D^3 u) \Big (\frac{\partial}{\partial \zeta},\frac{\partial}{\partial \bar{\zeta}},\frac{\partial}{\partial \zeta} \Big ) = O(|\zeta|^{\alpha\beta}).\] 
This completes the proof. \\

After these preparations, we now describe the estimate for the covariant derivative of $\partial \bar{\partial} u$. Following Yau \cite{Yau2}, we consider the function 
\[S = \sum_{\alpha,\beta,\gamma,\delta,\mu,\nu=1}^n \hat{g}^{\alpha\bar{\beta}} \, \hat{g}^{\gamma\bar{\delta}} \, \hat{g}^{\mu\bar{\nu}} \, (D^3 u) \Big ( \frac{\partial}{\partial z_\alpha},\frac{\partial}{\partial \bar{z}_\delta},\frac{\partial}{\partial z_\mu} \Big ) \, (D^3 u) \Big ( \frac{\partial}{\partial \bar{z}_\beta},\frac{\partial}{\partial z_\gamma},\frac{\partial}{\partial \bar{z}_\nu} \Big ).\] 
Here, the tensor $D^3 u$ represents the covariant derivative, with respect to $\omega$, of the $(1,1)$-form $\partial \bar{\partial} u$.

\begin{proposition}
\label{S.Holder}
Fix a point $p \in \Sigma$, and let $(z_1,\hdots,z_{n-1},\zeta)$ be complex coordinates around $p$. Then 
\[S(z_1,\hdots,z_{n-1},\zeta) = A(z_1,\hdots,z_{n-1}) + O(|\zeta|^{\alpha\beta}),\] 
where $A(z_1,\hdots,z_{n-1})$ is a continuous function on $\Sigma$.
\end{proposition}

\textbf{Proof.} 
It follows from Proposition \ref{covariant.derivatives.a} and Proposition \ref{covariant.derivatives.b} that the function $S$ can be written as a sum $S = S^{(1)} + S^{(2)}$, where $S^{(1)} \in \mathcal{C}^{,\alpha,\beta}$ and $S^{(2)} = O(|\zeta|^{\alpha\beta})$. From this, the assertion follows. \\

\begin{proposition} 
\label{c3}
Let $u \in \mathcal{C}^{2,\alpha,\beta}$ be a solution of $(\star_t)$ for some $t \in [0,1]$. Then $S \leq C$ for some uniform constant $C$.
\end{proposition}

\textbf{Proof.} 
Let us fix $\varepsilon > 0$ sufficiently small. It follows from Lemma \ref{barrier} that 
\[\sqrt{-1} \, \partial\bar{\partial}(|s|_h^{2\varepsilon}) + \omega \geq \varepsilon^2 \, |s|_h^{2\varepsilon} \, \sqrt{-1} \, \partial \log(|s|_h^2) \wedge \bar{\partial} \log(|s|_h^2),\] 
hence 
\[\Delta_{\hat{\omega}}(|s|_h^{2\varepsilon}) + \text{\rm tr}_{\hat{\omega}}(\omega) \geq \varepsilon^2 \, |s|_h^{2\varepsilon} \, \big | d \log(|s|_h^2) \big |_{\hat{g}}^2.\] 
We have shown in Section \ref{C2.estimate} that $a_1 \, \omega \leq \hat{\omega} \leq a_2 \, \omega$ for uniform constants $a_1,a_2 > 0$. Thus, we conclude that 
\begin{equation} 
\label{a}
\Delta_{\hat{\omega}} (|s|_h^{2\varepsilon}) \geq \delta \, |s|_h^{2\varepsilon-2\beta} - K 
\end{equation}
for uniform constants $\delta$ and $K$. We next consider the function $S$ defined above. Using the formula on page 266 of \cite{Aubin2}, we obtain 
\begin{align*} 
\Delta_{\hat{\omega}} S 
&\geq -C \, |D\partial\bar{\partial} u|_g^2 \, |\partial\bar{\partial} F|_g - C \, |D\partial\bar{\partial} u|_g \, |D\partial\bar{\partial} F|_g \\ 
&- C \, |D\partial\bar{\partial} u|_g^2 \, |R|_g - C \, |D\partial\bar{\partial} u|_g \, |DR|_g, 
\end{align*} 
where $D$ denotes the Levi-Civita connection associated with $\omega$ (see also \cite{Yau2}, Appendix A). Using the identity $\text{\rm Ric}_\omega = \sqrt{-1} \, \partial\bar{\partial} F$, we obtain 
\[\Delta_{\hat{\omega}} S \geq -C \, |D\partial\bar{\partial} u|_g^2 \, |R|_g - C \, |D\partial\bar{\partial} u|_g \, |DR|_g,\] 
hence 
\[\Delta_{\hat{\omega}} S \geq -C \, S \, |R|_g - C \, S^{\frac{1}{2}} \, |DR|_g.\] 
Moreover, it follows from Proposition \ref{curvature} and Proposition \ref{derivative.of.curvature} that $|R|_g \leq C$ and $|DR|_g \leq C \, |s|_h^{\varepsilon-\beta}$. Putting these facts together, we obtain 
\[\Delta_{\hat{\omega}} S \geq -C \, |s|_h^{\varepsilon-\beta} \, S^{\frac{1}{2}} - C \, S.\] 
Therefore, we can find a uniform constant $L$ such that 
\begin{equation}
\label{b}
\Delta_{\hat{\omega}} S \geq -\delta \, |s|_h^{2\varepsilon-2\beta} - L \, S, 
\end{equation} 
where $\delta$ is the constant appearing in (\ref{a}). Furthermore, equation (2.7) in \cite{Yau2} implies 
\begin{equation} 
\label{c}
\Delta_{\hat{\omega}} (\Delta_\omega u) \geq \kappa \, S - N 
\end{equation}
for uniform constants $\kappa$ and $N$. We next consider the function 
\[Q = S + \frac{L+1}{\kappa} \, \Delta_\omega u + |s|_h^{2\varepsilon}.\] 
Using the inequalities (\ref{a}), (\ref{b}), and (\ref{c}), we obtain 
\[\Delta_{\hat{\omega}} Q \geq S - \frac{L+1}{\kappa} \, N - K.\]
Since $\Delta_\omega u$ is uniformly bounded, we conclude that 
\begin{equation} 
\label{d}
\Delta_{\hat{\omega}} Q \geq Q - C 
\end{equation} 
for some uniform constant $C$. It follows from Proposition \ref{S.Holder} that the function $Q$ attains its maximum at some point in $M \setminus \Sigma$. Using (\ref{d}), we conclude that $Q$ is bounded from above by some uniform constant. From this, the assertion follows. \\

\begin{corollary}
\label{holder.bound}
Fix a point $p \in \Sigma$, and let $(z_1,\hdots,z_{n-1},\zeta)$ be complex coordinates around $p$. Then the functions 
\[\frac{\partial^2 u}{\partial z_k \, \partial \bar{z}_l}, \quad |\zeta|^{1-\beta} \, \frac{\partial^2 u}{\partial \zeta \, \partial \bar{z}_l}, \quad |\zeta|^{2-2\beta} \, \frac{\partial^2 u}{\partial \zeta \, \partial \bar{\zeta}}\] 
are Lipschitz continuous with respect to $\omega$, and the Lipschitz constants are bounded from above by some uniform constant $L$. In particular, $\|u\|_{\mathcal{C}^{2,\alpha,\beta}} \leq C$.
\end{corollary}

\textbf{Proof.} 
For abbreviation, let 
\[\varphi_{kl} = \frac{\partial^2 u}{\partial z_k \, \partial \bar{z}_l}, \quad \psi_l = |\zeta|^{1-\beta} \, \frac{\partial^2 u}{\partial \zeta \, \partial \bar{z}_l}, \quad \chi = |\zeta|^{2-2\beta} \, \frac{\partial^2 u}{\partial \zeta \, \partial \bar{\zeta}}.\] 
It follows from Proposition \ref{c3} that the covariant derivative of $\partial\bar{\partial} u$ is uniformly bounded, i.e. $|D\partial\bar{\partial} u|_g \leq C$ for some constant $C$. Moreover, we have 
\[\Big | D\frac{\partial}{\partial z_k} \Big |_g \leq C, \quad \Big | D \Big ( \zeta^{1-\beta} \, \frac{\partial}{\partial \zeta} \Big ) \Big |_g \leq C\] 
by Lemma \ref{christoffel}. Here, $D$ denotes the Levi-Civita connection associated with the background edge metric $\omega$. Note that the vector field $\zeta^{1-\beta} \, \frac{\partial}{\partial \zeta}$ is defined only locally. Since 
\[\varphi_{kl} = \partial\bar{\partial} u \Big ( \frac{\partial}{\partial z_k},\frac{\partial}{\partial \bar{z}_l} \Big ),\] 
we obtain 
\begin{equation} 
\label{varphi}
|d\varphi_{kl}|_g \leq C. 
\end{equation} 
Moreover, since 
\[\chi = \partial\bar{\partial} u \Big ( \zeta^{1-\beta} \, \frac{\partial}{\partial \zeta},\overline{\zeta^{1-\beta} \, \frac{\partial}{\partial \zeta}} \Big ),\] 
we obtain 
\begin{equation} 
\label{chi}
|d\chi|_g \leq C. 
\end{equation} 
Finally, using the identity 
\[\frac{\zeta^{1-\beta}}{|\zeta|^{1-\beta}} \, \psi_l = \partial\bar{\partial} u \Big ( \zeta^{1-\beta} \, \frac{\partial}{\partial \zeta},\frac{\partial}{\partial \bar{z}_l} \Big ),\] 
we obtain 
\begin{equation} 
\label{psi} 
\Big | d \Big ( \frac{\zeta^{1-\beta}}{|\zeta|^{1-\beta}} \, \psi_l \Big ) \Big |_g \leq C. 
\end{equation}
It follows from (\ref{varphi}) and (\ref{chi}) that the functions $\varphi_{kl}$ and $\chi$ are Lipschitz continuous, with Lipschitz constant at most $C$. It remains to bound the Lipschitz constant of $\psi_l$. Note that the function $\psi_l$ vanishes along $\Sigma$. Integrating the inequality (\ref{psi}) along radial line segments yields $|\psi_l(z_1,\hdots,z_{n-1},\zeta)| \leq C \, |\zeta|^\beta$. Using this inequality and (\ref{psi}), we obtain 
\begin{align*} 
|d\psi_l|_g 
&\leq \Big | d \Big ( \frac{\zeta^{1-\beta}}{|\zeta|^{1-\beta}} \, \psi_l \Big ) \Big |_g + |\psi_l| \, \Big | d \Big ( \frac{\zeta^{1-\beta}}{|\zeta|^{1-\beta}} \Big ) \Big |_g \\ 
&\leq C + C \, |\zeta|^\beta \, \Big | d \Big ( \frac{\zeta^{1-\beta}}{|\zeta|^{1-\beta}} \Big ) \Big |_g \\ 
&\leq C + C \, |\zeta|^{\beta-1} \, \sqrt{g^{\zeta\bar{\zeta}}}. 
\end{align*}
Consequently, $|d\psi_l|_g \leq C$ for some uniform constant $C$. This completes the proof. \\

\section{Proof of Theorem \ref{main.theorem}}

\label{proof.of.main.theorem}

We now finish the proof of Theorem \ref{main.theorem}. As above, we consider the equation 
\begin{equation} 
\tag{$\star_t$} (\omega + \sqrt{-1} \, \partial\bar{\partial} u)^n = e^{tF-c} \, \omega^n, 
\end{equation}
where $c$ is a constant. Moreover, let 
\[I = \{t \in [0,1]: \text{\rm $(\star_t)$ has a solution $(u,c) \in \mathcal{C}^{2,\alpha,\beta} \times \mathbb{R}$}\}.\] 
The equation $(\star_0)$ has the trivial solution $u = 0$; in particular, $I$ is non-empty. It follows from Corollary \ref{holder.bound} that any solution of $(\star_t)$ is uniformly bounded in $\mathcal{C}^{2,\alpha,\beta}$. Consequently, $I$ is closed. We next show that $I$ is open.

\begin{proposition}
The set $I$ is an open subset of $[0,1]$.
\end{proposition} 

\textbf{Proof.} 
Fix a point $t_0 \in I$. Moreover, let $u$ be a solution of $(\star_{t_0})$, and let $\hat{\omega} = \omega + \sqrt{-1} \, \partial\bar{\partial} u$. By Theorem \ref{fredholm}, the operator $\Delta_{\hat{\omega}}: \mathcal{C}^{2,\alpha,\beta} \to \mathcal{C}^{,\alpha,\beta}$ is Fredholm with Fredholm index zero. Consequently, the operator $\Delta_{\hat{\omega}}$ has one-dimensional cokernel and its range is given by $\big \{ f \in \mathcal{C}^{,\alpha,\beta}: \int_M f \, \hat{\omega}^n = 0 \big \}$. Hence, the assertion follows from the implicit function theorem. \\

\begin{corollary}
\label{existence}
We have $I = [0,1]$. In particular, the equation $(\star_1)$ admits a solution $u \in \mathcal{C}^{2,\alpha,\beta}$.
\end{corollary} 

In the remainder of this section, we assume that $u \in \mathcal{C}^{2,\alpha,\beta}$ is a solution of $(\star_1)$. Then the associated K\"ahler metric $\hat{\omega} = \omega + \sqrt{-1} \, \partial\bar{\partial} u$ is Ricci flat. It remains to show that $\hat{\omega}$ has bounded curvature. The proof of this statement can be broken down into several lemmata.

\begin{lemma} 
\label{fourth.order.derivatives.a}
Let $u \in \mathcal{C}^{2,\alpha,\beta}$ be a solution of $(\star_1)$. Then we have 
\begin{align*} 
&\frac{\partial^2}{\partial z_i \, \partial \bar{z}_j} \hat{g}_{k\bar{l}} \in \mathcal{C}^{,\alpha,\beta}, \\ 
&|\zeta|^{1-\beta} \, \frac{\partial}{\partial z_i \, \partial \bar{\zeta}} \hat{g}_{k\bar{l}} \in \mathcal{C}_0^{,\alpha,\beta}, \\ 
&|\zeta|^{1-\beta} \, \frac{\partial}{\partial \zeta \, \partial \bar{z}_j} \hat{g}_{k\bar{l}} \in \mathcal{C}_0^{,\alpha,\beta}, \\ 
&|\zeta|^{2-2\beta} \, \frac{\partial}{\partial \zeta \, \partial \bar{\zeta}} \hat{g}_{k\bar{l}} \in \mathcal{C}^{,\alpha,\beta} 
\end{align*} 
for all $i,j,k,l \in \{1,\hdots,n-1\}$.
\end{lemma}

\textbf{Proof.} 
Let us fix two integers $k,l \in \{1,\hdots,n-1\}$. Since $u$ is a solution of $(\star_1)$, the Ricci tensor of $\hat{g}$ vanishes. This implies 
\[\Delta_{\hat{\omega}} \hat{g}_{k\bar{l}} - \sum_{\alpha,\beta,\mu,\nu=1}^n \hat{g}^{\alpha\bar{\beta}} \, \hat{g}^{\mu\bar{\nu}} \, \frac{\partial}{\partial z_k} \hat{g}_{\alpha\bar{\nu}} \, \frac{\partial}{\partial \bar{z}_l} \hat{g}_{\mu\bar{\beta}} = 0.\] 
By Lemma \ref{third.order.derivatives.a}, the $(1,1)$-forms $\frac{\partial}{\partial z_k} \hat{\omega}$ and $\frac{\partial}{\partial \bar{z}_l} \hat{\omega}$ are of class $\mathcal{C}^{,\alpha,\beta}$. Consequently, we can find a function $f \in \mathcal{C}^{,\alpha,\beta}$ such that $\Delta_{\hat{\omega}} \hat{g}_{k\bar{l}} = f$ away from $\Sigma$. Hence, the assertion follows from Corollary \ref{regularity.3}. \\

\begin{lemma} 
\label{fourth.order.derivatives.b}
Let $u \in \mathcal{C}^{2,\alpha,\beta}$ be a solution of $(\star_1)$. Then we have 
\[|\zeta|^{2-2\beta} \, \Big ( \frac{\partial^4 u}{\partial \zeta \, \partial \bar{z}_j \, \partial \zeta \, \partial \bar{z}_l} + \frac{1-\beta}{\zeta} \, \frac{\partial^3 u}{\partial \bar{z}_j \, \partial \zeta \, \partial \bar{z}_l} \Big ) = O(|\zeta|^{\alpha\beta})\] 
for all $j,l \in \{1,\hdots,n-1\}$.
\end{lemma}

\textbf{Proof.} 
Fix two integers $j,l \in \{1,\hdots,n-1\}$. By assumption, we have $u \in \mathcal{C}^{2,\alpha,\beta}$. In particular, we have $\Delta_\Omega u \in \mathcal{C}^{,\alpha,\beta}$. Using the results in \cite{Donaldson}, we conclude that $\frac{\partial^2 u}{\partial \bar{z}_j \, \partial \bar{z}_l} \in \mathcal{C}^{,\alpha,\beta}$. Since $u$ is a solution of $(\star_1)$, we have 
\[\frac{\partial^2}{\partial \bar{z}_j \, \partial \bar{z}_l} \log \hat{\omega}^n = \frac{\partial^2}{\partial \bar{z}_j \, \partial \bar{z}_l} (\log \omega^n + F).\] 
This implies 
\begin{align*} 
&\Delta_{\hat{\omega}} \Big ( \frac{\partial^2 u}{\partial \bar{z}_j \, \partial \bar{z}_l} \Big ) + \sum_{\mu,\nu=1}^n \hat{g}^{\mu\bar{\nu}} \, \frac{\partial^2}{\partial \bar{z}_j \, \partial \bar{z}_l} g_{\mu\bar{\nu}} - \sum_{\alpha,\beta,\mu,\nu=1}^n \hat{g}^{\alpha\bar{\beta}} \, \hat{g}^{\mu\bar{\nu}} \, \frac{\partial}{\partial \bar{z}_j} \hat{g}_{\alpha\bar{\nu}} \, \frac{\partial}{\partial \bar{z}_l} \hat{g}_{\mu\bar{\beta}} \\ 
&= \sum_{\mu,\nu=1}^n \hat{g}^{\mu\bar{\nu}} \, \frac{\partial^2}{\partial \bar{z}_j \, \partial \bar{z}_l} \hat{g}_{\mu\bar{\nu}} - \sum_{\alpha,\beta,\mu,\nu=1}^n \hat{g}^{\alpha\bar{\beta}} \, \hat{g}^{\mu\bar{\nu}} \, \frac{\partial}{\partial \bar{z}_j} \hat{g}_{\alpha\bar{\nu}} \, \frac{\partial}{\partial \bar{z}_l} \hat{g}_{\mu\bar{\beta}} \\ 
&= \frac{\partial^2}{\partial \bar{z}_j \, \partial \bar{z}_l} (\log \omega^n + F). 
\end{align*} 
By Lemma \ref{third.order.derivatives.a}, the $(1,1)$-forms $\frac{\partial}{\partial \bar{z}_j} \hat{\omega}$ and $\frac{\partial}{\partial \bar{z}_l} \hat{\omega}$ are of class $\mathcal{C}^{,\alpha,\beta}$. Moreover, it is easy to see that the function $\frac{\partial^2}{\partial \bar{z}_j \, \partial \bar{z}_l} (\log \omega^n + F)$ is of class $\mathcal{C}^{,\alpha,\beta}$. Consequently, there exists a function $f \in \mathcal{C}^{,\alpha,\beta}$ such that $\Delta_{\hat{\omega}} \big ( \frac{\partial^2 u}{\partial \bar{z}_j \, \partial \bar{z}_l} \big ) = f$ away from $\Sigma$. Using Corollary \ref{regularity.3}, we conclude that 
\[\frac{\partial^2 u}{\partial \bar{z}_j \, \partial \bar{z}_l} \in \mathcal{C}^{2,\alpha,\beta}.\] 
In particular, we have 
\[|\zeta|^{2-2\beta} \, \frac{\partial^4 u}{\partial \zeta \, \partial \bar{\zeta} \, \partial \bar{z}_j \, \partial \bar{z}_l} \in \mathcal{C}^{,\alpha,\beta}.\] 
Hence, the assertion follows from Proposition \ref{auxiliary.result}. \\

\begin{proposition}
\label{bounded.curvature}
Let $u \in \mathcal{C}^{2,\alpha,\beta}$ be a solution of $(\star_1)$, and let $\hat{R}$ denote the Riemann curvature tensor of $\hat{g}$. Then $\sup_{M \setminus \Sigma} |\hat{R}|_{\hat{g}} < \infty$.
\end{proposition}

\textbf{Proof.} As in Section \ref{background.metric}, we consider the map 
\[\Psi: (z_1,\hdots,z_{n-1},z_n) \mapsto (z_1,\hdots,z_{n-1},z_n^{\frac{1}{\beta}}).\] 
We have shown in Section \ref{C3.estimate} that $|D\hat{\omega}|_g = |D\partial\bar{\partial} u|_g = O(1)$. This implies 
\[\sum_{\alpha,\gamma,\delta=1}^n \Big | \frac{\partial}{\partial z_\alpha} \Psi^* \hat{g} \Big ( \frac{\partial}{\partial z_\gamma},\frac{\partial}{\partial \bar{z}_\delta} \Big ) \Big | = O(1).\] 
Using Lemma \ref{fourth.order.derivatives.a}, we obtain 
\[\sum_{i,j,k,l=1}^{n-1} \Big | \frac{\partial^2}{\partial z_i \, \partial \bar{z}_j} \Psi^* \hat{g} \Big ( \frac{\partial}{\partial z_k},\frac{\partial}{\partial \bar{z}_l} \Big ) \Big | = O(1)\] 
and 
\[\sum_{j,k,l=1}^{n-1} \Big | \frac{\partial^2}{\partial z_n \, \partial \bar{z}_j} \Psi^* \hat{g} \Big ( \frac{\partial}{\partial z_k},\frac{\partial}{\partial \bar{z}_l} \Big ) \Big | = O(1).\] 
Consequently, the Riemann curvature tensor of $\hat{g}$ satisfies 
\begin{equation} 
\label{riem.1}
\sum_{i,j,k,l=1}^{n-1} \Big | \hat{R} \Big ( \frac{\partial}{\partial z_i},\frac{\partial}{\partial \bar{z}_j},\frac{\partial}{\partial z_k},\frac{\partial}{\partial \bar{z}_l} \Big ) \Big | = O(1)
\end{equation} 
and 
\begin{equation} 
\label{riem.2} 
\sum_{j,k,l=1}^{n-1} \Big | \hat{R} \Big ( \frac{\partial}{\partial z_n},\frac{\partial}{\partial \bar{z}_j},\frac{\partial}{\partial z_k},\frac{\partial}{\partial \bar{z}_l} \Big ) \Big | = O(1). 
\end{equation}
Moreover, Lemma \ref{fourth.order.derivatives.b} implies that 
\[\sum_{j,l=1}^{n-1} \Big | \frac{\partial^4}{\partial z_n \, \partial \bar{z}_j \, \partial z_n \, \partial \bar{z}_l} (u \circ \Psi) \Big | = O(1).\] 
Using this inequality and the estimate 
\[\sum_{j,l=1}^{n-1} \Big | \frac{\partial^2}{\partial z_n \, \partial \bar{z}_j} \Psi^* g \Big ( \frac{\partial}{\partial z_n},\frac{\partial}{\partial \bar{z}_l} \Big ) \Big | = O(1),\] 
we obtain 
\[\sum_{j,l=1}^{n-1} \Big | \frac{\partial^2}{\partial z_n \, \partial \bar{z}_j} \Psi^* \hat{g} \Big ( \frac{\partial}{\partial z_n},\frac{\partial}{\partial \bar{z}_l} \Big ) \Big | = O(1).\] 
Thus, we conclude that 
\begin{equation} 
\label{riem.3} 
\sum_{j,l=1}^{n-1} \Big | \hat{R} \Big ( \frac{\partial}{\partial z_n},\frac{\partial}{\partial \bar{z}_j},\frac{\partial}{\partial z_n},\frac{\partial}{\partial \bar{z}_l} \Big ) \Big | = O(1). 
\end{equation} 
Since the Ricci tensor of $\hat{g}$ vanishes, the inequalities (\ref{riem.1}), (\ref{riem.2}), and (\ref{riem.3}) imply that 
\[\sum_{\alpha,\beta,\gamma,\delta=1}^n \Big | \hat{R} \Big ( \frac{\partial}{\partial z_\alpha},\frac{\partial}{\partial \bar{z}_\beta},\frac{\partial}{\partial z_\gamma},\frac{\partial}{\partial \bar{z}_\delta} \Big ) \Big | = O(1).\] 
Since the metric $\hat{g}$ is uniformly equivalent to $g$, we conclude that $|\hat{R}|_{\hat{g}} = O(1)$. This completes the proof of Proposition \ref{bounded.curvature}. 

\appendix 

\section{An auxiliary result}

\begin{proposition} 
\label{auxiliary.result}
Fix real numbers $\alpha \in (0,1)$ and $\beta \in (0,\frac{1}{2})$ such that $\alpha\beta < 1-2\beta$. Assume that $\tilde{f}$ is a function defined on the unit disk $B_1(0) \subset \mathbb{C}$ of class $C^\alpha$. Moreover, suppose that $v$ is a smooth function defined on the punctured disk $B_1(0) \setminus \{0\} \subset \mathbb{C}$ satisfying $\sup_{B_1(0) \setminus \{0\}} |v| < \infty$ and 
\[|\zeta|^{2-2\beta} \, \frac{\partial^2 v}{\partial \zeta \, \partial \bar{\zeta}} = \tilde{f}(|\zeta|^{\beta-1} \, \zeta)\] 
on $B_1(0) \setminus \{0\}$. Then 
\[|\zeta|^{2-2\beta} \, \Big ( \frac{\partial^2 v}{\partial \zeta^2} + \frac{1-\beta}{\zeta} \, \frac{\partial v}{\partial \zeta} \Big ) = O(|\zeta|^{\alpha\beta}).\] 
\end{proposition}

\textbf{Proof.} Let 
\[w(\zeta) = v(\zeta) - \frac{\tilde{f}(0)}{\beta^2} \, |\zeta|^{2\beta}\] 
and 
\[h(\zeta) = |\zeta|^{2\beta-2} \, (\tilde{f}(|\zeta|^{\beta-1} \, \zeta) - \tilde{f}(0)).\] 
Since $\tilde{f}$ is of class $C^\alpha$, we have $|h(\zeta)| \leq C \, |\zeta|^{2\beta-2+\alpha\beta}$. Moreover, we have 
\[\frac{\partial^2 w}{\partial \zeta \, \partial \bar{\zeta}} = h\] 
on $B_1(0) \setminus \{0\}$. Since $\sup_{B_1(0) \setminus \{0\}} |w| < \infty$, the previous equation holds in the distributional sense. Using Green's formula, we obtain 
\begin{align*} 
\Big | \frac{\partial w}{\partial \zeta}(\zeta_0) \Big | 
&\leq C \int_{B_{\frac{1}{2}}(0)} |\zeta-\zeta_0|^{-1} \, |h(\zeta)| \, d\zeta + C \\ 
&\leq C \int_{B_{\frac{1}{2}}(0)} |\zeta-\zeta_0|^{-1} \, |\zeta|^{2\beta-2+\alpha\beta} \, d\zeta + C \\ 
&\leq C \, |\zeta_0|^{2\beta-1+\alpha\beta}. 
\end{align*} 
for $\zeta_0 \in B_{1/4}(0)$. Consequently, 
\[|\zeta|^{1-2\beta} \, \frac{\partial w}{\partial \zeta} = O(|\zeta|^{\alpha\beta}).\] 
We now consider a point $\zeta_0 \in B_{1/4}(0)$. Since $\tilde{f}$ is of class $C^\alpha$, we have 
\begin{align*} 
|h(\zeta_1) - h(\zeta_2)| 
&\leq |\zeta_1|^{2\beta-2} \, \big | \tilde{f}(|\zeta_1|^{\beta-1} \, \zeta_1) - \tilde{f}(|\zeta_2|^{\beta-1} \, \zeta_2) \big | \\ 
&+ \big | |\zeta_1|^{2\beta-2} - |\zeta_2|^{2\beta-2} \big | \, \big | \tilde{f}(|\zeta_2|^{\beta-1} \, \zeta_2) - \tilde{f}(0) \big | \\ 
&\leq C \, |\zeta_1|^{2\beta-2} \, \big | |\zeta_1|^{\beta-1} \, \zeta_1 - |\zeta_2|^{\beta-1} \, \zeta_2|^\alpha \\ 
&+ C \, \big | |\zeta_1|^{2\beta-2} - |\zeta_2|^{2\beta-2} \big | \, |\zeta_2|^{\alpha\beta}
\end{align*} 
for all $\zeta_1,\zeta_2 \in B_{|\zeta_0|/2}(\zeta_0)$. This implies 
\begin{align*} 
|h(\zeta_1) - h(\zeta_2)| 
&\leq C \, |\zeta_0|^{2\beta-2+\alpha(\beta-1)} \, |\zeta_1-\zeta_2|^\alpha \\ 
&+ C \, |\zeta_0|^{2\beta-3+\alpha\beta} \, |\zeta_1-\zeta_2| 
\end{align*} 
for all $\zeta_1,\zeta_2 \in B_{|\zeta_0|/2}(\zeta_0)$. Therefore, 
\[[h]_{C^\alpha(B_{|\zeta_0|/2}(\zeta_0))} \leq C \, |\zeta_0|^{2\beta-2+\alpha(\beta-1)}.\] 
Using standard interior estimates, we obtain 
\begin{align*} 
\Big | \frac{\partial^2 w}{\partial \zeta^2}(\zeta_0) \Big | 
&\leq C \, |\zeta_0|^{-1} \, \sup_{B_{|\zeta_0|/2}(\zeta_0)} \Big | \frac{\partial w}{\partial \zeta} \Big | \\ 
&+ C \, \sup_{B_{|\zeta_0|/2}(\zeta_0)} |h| + C \, |\zeta_0|^\alpha \, [h]_{C^\alpha(B_{|\zeta_0|/2}(\zeta_0))}. 
\end{align*}
Putting these facts together, we conclude that 
\[\Big | \frac{\partial^2 w}{\partial \zeta^2}(\zeta_0) \Big | \leq C \, |\zeta_0|^{2\beta-2+\alpha\beta}.\] 
Consequently, 
\[|\zeta|^{2-2\beta} \, \frac{\partial^2 w}{\partial \zeta^2} = O(|\zeta|^{\alpha\beta}).\] 
Putting these facts together, we obtain  
\[|\zeta|^{2-2\beta} \, \Big ( \frac{\partial^2 v}{\partial \zeta^2} + \frac{1-\beta}{\zeta} \, \frac{\partial v}{\partial \zeta} \Big ) = |\zeta|^{2-2\beta} \, \Big ( \frac{\partial^2 w}{\partial \zeta^2} + \frac{1-\beta}{\zeta} \, \frac{\partial w}{\partial \zeta} \Big ) = O(|\zeta|^{\alpha\beta}).\]


\begin{thebibliography}{99}
\bibitem{Aubin1}
T.~Aubin, \textit{\'Equations du type Monge-Amp\`ere sur les vari\'et\'es K\"ahleriennes compactes,} C.R. Acad. Sci. Paris S\'er. A-B 283, A119--A121 (1976)

\bibitem{Aubin2}
T.~Aubin, \textit{Some nonlinear problems in Riemannian geometry,} Springer Monographs in Mathematics, Springer, 1998

\bibitem{Donaldson}
S.~Donaldson, \textit{K\"ahler metrics with cone singularities along a divisor,} arxiv:1102.1196

\bibitem{Gilbarg-Trudinger}
D.~Gilbarg and N.~Trudinger, \textit{Elliptic Partial Differential Equations of Second Order,} Springer-Verlag, 2001

\bibitem{Jeffres1}
T.~Jeffres, \textit{Schwarz lemma for K\"ahler cone metrics,} Internat. Math. Res. Notices 371--382 (2000)

\bibitem{Jeffres2}
T.~Jeffres, \textit{Uniqueness of K\"ahler-Einstein cone metrics,} Publ. Mat. 44, 437--448 (2000)

\bibitem{Mazzeo1}
R.~Mazzeo, \textit{Elliptic theory of differential edge operators, I,} Comm. Partial Differential Equations 16, 1615--1664 (1991)

\bibitem{Mazzeo2}
R.~Mazzeo, \textit{K\"ahler-Einstein metrics singular along a smooth divisor,} Journ\'ees \'Equations aux D\'eriv\'ees Partielles (Saint-Jean-de-Monts, 1999), Exp. No. VI, 1--10

\bibitem{Simon}
L.~Simon, \textit{Schauder estimates by scaling,} Calc. Var. PDE 5, 391--407 (1997)

\bibitem{Tian}
G.~Tian, \textit{Canonical metrics in K\"ahler geometry,} Birkh\"auser, 2000

\bibitem{Tian-Yau1}
G.~Tian and S.T.~Yau, \textit{Complete K\"ahler manifolds with zero Ricci curvature I,} J. Amer. Math. Soc. 3, 579--609 (1990)

\bibitem{Tian-Yau2}
G.~Tian and S.T.~Yau, \textit{Complete K\"ahler manifolds with zero Ricci curvature II,} Invent. Math. 106, 27--60 (1991)

\bibitem{Troyanov}
M.~Troyanov, \textit{Prescribing curvature on compact surfaces with conical singularities,} Trans. Amer. Math. Soc. 325, 793--821 (1991)

\bibitem{Yau1}
S.T.~Yau, \textit{On Calabi's conjecture and some new results in algebraic geometry,} Proc. Natl. Acad. Sci. USA 74, 1798--1799 (1977)

\bibitem{Yau2}
S.T.~Yau, \textit{On the Ricci curvature of a compact K\"ahler manifold and the complex Monge-Amp\`ere equation,} Comm. Pure Appl. Math. 31, 339--411 (1978)
\end{thebibliography}
\end{document}